\documentclass{article}

\usepackage{amssymb,latexsym,amsmath}
\usepackage[pdftex]{graphicx}

\usepackage{hyperref}

\begin{document}

\pagestyle{plain}

\newtheorem{theorem}{Theorem}[section]

\newtheorem{proposition}[theorem]{Proposition}

\newtheorem{lema}[theorem]{Lemma}

\newtheorem{corollary}[theorem]{Corollary}

\newtheorem{definition}[theorem]{Definition}

\newtheorem{remark}[theorem]{Remark}

\newtheorem{exempl}{Example}[section]

\newenvironment{exemplu}{\begin{exempl}  \em}{\hfill $\square$

\end{exempl}}  \vspace{.5cm}

\newcommand{\ea}{\mbox{{\bf a}}  \vspace{.5cm}}

\newcommand{\eu}{\mbox{{\bf u}}  \vspace{.5cm}}

\newcommand{\ueu}{\underline{\eu}}  \vspace{.5cm}

\newcommand{\ueo}{\overline{u}}  \vspace{.5cm}

\newcommand{\oeu}{\overline{\eu}}  \vspace{.5cm}

\newcommand{\ew}{\mbox{{\bf w}}  \vspace{.5cm}}

\newcommand{\ef}{\mbox{{\bf f}}  \vspace{.5cm}}

\newcommand{\eF}{\mbox{{\bf F}}  \vspace{.5cm}}

\newcommand{\eC}{\mbox{{\bf C}}  \vspace{.5cm}}

\newcommand{\en}{\mbox{{\bf n}}  \vspace{.5cm}}

\newcommand{\eT}{\mbox{{\bf T}}  \vspace{.5cm}}

\newcommand{\eL}{\mbox{{\bf L}}  \vspace{.5cm}}

\newcommand{\eR}{\mbox{{\bf R}}  \vspace{.5cm}}

\newcommand{\eV}{\mbox{{\bf V}}  \vspace{.5cm}}

\newcommand{\eU}{\mbox{{\bf U}}  \vspace{.5cm}}

\newcommand{\ev}{\mbox{{\bf v}}  \vspace{.5cm}}

\newcommand{\eve}{\mbox{{\bf e}}  \vspace{.5cm}}

\newcommand{\uev}{\underline{\ev}}  \vspace{.5cm}

\newcommand{\eY}{\mbox{{\bf Y}}  \vspace{.5cm}}

\newcommand{\eK}{\mbox{{\bf K}}  \vspace{.5cm}}

\newcommand{\eP}{\mbox{{\bf P}}  \vspace{.5cm}}

\newcommand{\eS}{\mbox{{\bf S}}  \vspace{.5cm}}

\newcommand{\eJ}{\mbox{{\bf J}}  \vspace{.5cm}}

\newcommand{\eB}{\mbox{{\bf B}}  \vspace{.5cm}}

\newcommand{\eH}{\mbox{{\bf H}}  \vspace{.5cm}}

\newcommand{\leb}{\mathcal{ L}^{n}}  \vspace{.5cm}

\newcommand{\eI}{\mathcal{ I}}  \vspace{.5cm}

\newcommand{\eE}{\mathcal{ E}}  \vspace{.5cm}

\newcommand{\hen}{\mathcal{H}^{n-1}}  \vspace{.5cm}

\newcommand{\eBV}{\mbox{{\bf BV}}  \vspace{.5cm}}

\newcommand{\eA}{\mbox{{\bf A}}  \vspace{.5cm}}

\newcommand{\eSBV}{\mbox{{\bf SBV}}  \vspace{.5cm}}

\newcommand{\eBD}{\mbox{{\bf BD}}  \vspace{.5cm}}

\newcommand{\eSBD}{\mbox{{\bf SBD}}  \vspace{.5cm}}

\newcommand{\ecs}{\mbox{{\bf X}}  \vspace{.5cm}}

\newcommand{\eg}{\mbox{{\bf g}}  \vspace{.5cm}}

\newcommand{\paromega}{\partial \Omega}

\newcommand{\gau}{\Gamma_{u}}  \vspace{.5cm}

\newcommand{\gaf}{\Gamma_{f}}  \vspace{.5cm}

\newcommand{\sig}{{\bf \sigma}}  \vspace{.5cm}

\newcommand{\gac}{\Gamma_{\mbox{{\bf c}}  \vspace{.5cm}}}  \vspace{.5cm}

\newcommand{\deu}{\dot{\eu}}  \vspace{.5cm}

\newcommand{\dueu}{\underline{\deu}}  \vspace{.5cm}

\newcommand{\dev}{\dot{\ev}}  \vspace{.5cm}

\newcommand{\duev}{\underline{\dev}}  \vspace{.5cm}

\newcommand{\weak}{\stackrel{w}{\approx}}  \vspace{.5cm}

\newcommand{\mild}{\stackrel{m}{\approx}}  \vspace{.5cm}

\newcommand{\lrightarrow}{\stackrel{L}{\rightarrow}}  \vspace{.5cm}

\newcommand{\rrightarrow}{\stackrel{R}{\rightarrow}}  \vspace{.5cm}

\newcommand{\strong}{\stackrel{s}{\approx}}  \vspace{.5cm}

\newcommand{\weakdown}{\rightharpoondown}

\newcommand{\opg}{\stackrel{\mathfrak{g}}  \vspace{.5cm}{\cdot}}  \vspace{.5cm}

\newcommand{\opunu}{\stackrel{1}{\cdot}}  \vspace{.5cm}
\newcommand{\opdoi}{\stackrel{2}{\cdot}}  \vspace{.5cm}

\newcommand{\opn}{\stackrel{\mathfrak{n}}  \vspace{.5cm}{\cdot}}  \vspace{.5cm}
\newcommand{\opx}{\stackrel{x}{\cdot}}  \vspace{.5cm}

\newcommand{\tr}{\ \mbox{tr}}  \vspace{.5cm}

\newcommand{\Ad}{\ \mbox{Ad}}  \vspace{.5cm}

\newcommand{\ad}{\ \mbox{ad}}  \vspace{.5cm}

\renewcommand{\contentsname}{ }

\title{On graphic lambda calculus and the dual of the graphic beta move}

\author{Marius Buliga \\ 
\\
Institute of Mathematics, Romanian Academy \\
P.O. BOX 1-764, RO 014700\\
Bucure\c sti, Romania\\
{\footnotesize Marius.Buliga@imar.ro}}  \vspace{.5cm}

\date{This version: 04.02.2013}

\maketitle

\begin{abstract}
This is a short description of graphic lambda calculus, with special emphasis on a duality suggested by the two different appearances of knot diagrams, in  lambda calculus  and emergent algebra sectors of the graphic lambda calculus respectively. This duality leads to the introduction of the dual of the graphic beta move. While the graphic beta move corresponds to beta reduction in untyped lambda calculus, the dual graphic beta move appears in relation to emergent algebras. 
\end{abstract}

\tableofcontents

\newpage

\paragraph{Sources, to read if necessary:}       

The following list of papers, with the format (link to arxiv, bibliography citation, title):
\begin{enumerate}   
\item[-]   \href{http://arxiv.org/abs/1211.1604}{arXiv:1211.1604},  \cite{graphiclambdaknots}   - Graphic lambda calculus and knot diagrams   
\item[-]       \href{http://arxiv.org/abs/1207.0332}{arXiv:1207.0332}, \cite{graphiclambdamoves} - Local and global moves on locally planar trivalent graphs, lambda calculus and lambda-Scale   
\item[-]       \href{http://arxiv.org/abs/1205.0139}{arXiv:1205.0139}, \cite{lambdascale} - Lambda-Scale, a lambda calculus for spaces with dilations   
\item[-] \href{http://arxiv.org/abs/1005.5031}{arXiv:1005.5031}, \cite{buligabraided} - Braided spaces with dilations and sub-riemannian symmetric spaces
\item[-]       \href{http://arxiv.org/abs/0907.1520}{arXiv:0907.1520}, \cite{buligairq}  - Emergent algebras.     
\end{enumerate}

\begin{enumerate}
\item[-]  The online tutorial \href{http://chorasimilarity.wordpress.com/graphic-lambda-calculus/}{Graphic lambda calculus} and the posts with the tag     \href{http://chorasimilarity.wordpress.com/tag/graphic-lambda-calculus/}{graphic lambda calculus}    from the     \href{http://chorasimilarity.wordpress.com/}{chorasimilarity}    blog,    
	\item[-]   The page    \href{http://imar.ro/~mbuliga/buliga_sim.html}{Emergent algebras / Computing with space / Graphical calculi}   
\end{enumerate}

\paragraph{What is graphic lambda calculus?}       

Graphic lambda calculus is a formalism working with oriented, locally planar, trivalent (or univalent) graphs, with decorated nodes. It has a number of moves  acting on such graphs, which can be local or global (in the sense of  \href{http://arxiv.org/abs/1207.0332}{arXiv:1207.0332}, \cite{graphiclambdamoves} section 2).

It contains  differential calculus in metric spaces (via the formalism of emergent algebras 
\cite{buligairq}, itself constructed over the geometry and analysis of metric spaces with dilations, introduced in \cite{buligadil1}, see especially the section 4, on the use of binary decorates trees), untyped lambda calculus (explained in \href{http://arxiv.org/abs/1207.0332}{arXiv:1207.0332}, \cite{graphiclambdamoves}) and that part of knot theory which can be expressed by using knot diagrams (explained in \href{http://arxiv.org/abs/1211.1604}{arXiv:1211.1604},  \cite{graphiclambdaknots}).

\paragraph{What is the motivation?}       

I wanted to have a lambda calculus formalism which does not need alpha renaming, even better, without  names for variables or terms.

Also, I wanted to understand graphic manipulations of knot diagrams related to emergent algebras  (explained in the paper   \href{http://arxiv.org/abs/1103.6007}{arXiv:1103.6007}, \cite{buligachora}, for a more algebraic approach see \href{http://arxiv.org/abs/1005.5031}{arXiv:1005.5031}, \cite{buligabraided})  as a kind of computation. This is achieved by the extended graphic beta move.

\paragraph{Are there other graphic lambda calculi?}

Yes, see for example the page   \href{http://dkeenan.com/Lambda/}{To dissect a mockingbird: a graphical notation for the lambda calculus with animated reduction} and go to the end of the page for links to  bibliographic information.

However, the graphic lambda calculus presented here is different, for several reasons, among them because of the connections with knot diagrams (first explained in the post \href{http://chorasimilarity.wordpress.com/2012/10/26/3d-crossings-in-graphic-lambda-calculus/}{3D grossings in graphic lambda calculus}, the beta reduction is a braiding operation), because of the attention to local and global moves and sets of graphs (terms) and because of the connection with emergent algebras and differential calculus in metric spaces.

A detailed discussion about the relations between this graphic lambda calculus and others would be very useful.

\paragraph{Sectors of graphic lambda calculus.} A sector of the graphic lambda calculus is a pair:
\begin{enumerate}
\item[-]  a set of graphs, defined by a local or global condition,   
\item[-]   a set of moves, either  from the list of available moves, or moves obtained as a concatenation of the ones listed here.   
\end{enumerate}

The name “graphic lambda calculus” comes from the fact that there it has untyped lambda calculus as a sector. In fact, there are three important sectors of graphic lambda calculus:
\begin{enumerate}
\item[-]     untyped lambda calculus sector, which contains all graphs in $GRAPH$ which are obtained from \href{http://en.wikipedia.org/wiki/Lambda_calculus}{untyped lambda calculus}      terms  by \href{http://chorasimilarity.wordpress.com/2012/12/21/conversion-of-lambda-calculus-terms-into-graphs/}{the algorithm described here}, \cite{graphiclambdamoves} section 3.  The moves of this sector are: \href{http://chorasimilarity.wordpress.com/2012/12/13/the-graphic-beta-move-with-details/}{graphic beta move} section \ref{graphicbetawithdetails}, \href{http://chorasimilarity.wordpress.com/2012/12/17/fan-out-moves-co-comm-co-assoc-global-fan-out-local-fan-out/}{fan-out moves} section \ref{fanoutmoves}, \href{http://chorasimilarity.wordpress.com/2012/12/17/pruning-moves/}{pruning moves} section \ref{pruningmoves} and elimination of loops. The article \href{http://arxiv.org/abs/1207.0332}{arXiv:1207.0332} \cite{graphiclambdamoves} describes this sector in detail.     
\item[-]     emergent algebra sector, which contain all graphs in $GRAPH$ described in the article \href{http://arxiv.org/abs/1103.6007}{arXiv:1103.6007} \cite{buligachora} , via the \href{http://chorasimilarity.wordpress.com/2012/10/31/3d-crossings-in-emergent-algebras/}{emergent algebra crossing}  macros, and the following moves: \href{http://chorasimilarity.wordpress.com/tag/extended-beta-move/}{dual graphic beta move} (which forms with the graphic move the extended beta move), \href{http://chorasimilarity.wordpress.com/2012/12/17/fan-out-moves-co-comm-co-assoc-global-fan-out-local-fan-out/}{fan-out moves}, \href{http://chorasimilarity.wordpress.com/2012/12/17/pruning-moves/}{pruning moves}, \href{http://chorasimilarity.wordpress.com/2012/12/20/emergent-algebra-moves-r1a-r1b-and-ext2/}{emergent algebra moves} section \ref{emergentalgebramoves}.     
\item[-]     knot and tangle diagrams sector, defined by using \href{http://chorasimilarity.wordpress.com/2012/09/06/graphic-beta-rule-as-braiding/}{crossings in lambda calculus}  macro      and   the Reidemeister moves as described in the post \href{http://chorasimilarity.wordpress.com/2012/11/05/generating-set-of-reidemeister-moves-for-graphic-lambda-crossings/}{Generating set of Reidemeister moves for graphic lambda crossings}, which are composite moves obtained from the graphic beta move and some local fan-out moves. The article \href{http://arxiv.org/abs/1211.1604}{arXiv:1211.1604} \cite{graphiclambdaknots} describes this sector in detail.     
\end{enumerate}
The separation of the set $GRAPH$ of graphs into sectors is artificial, distinctions between sectors blur under a more attentive look. For example, untyped lambda calculus can be seen as represented by the untyped lambda calculus sector, but one can use a sector which exploits graphical versions of  combinators in order to obtain a different representation (not explained in this paper). Likewise, the knot and tangle diagram sector is constructed over the "lambda calculus crossing macro" described in section 
\ref{graphicbetaruleasbraiding}. One can use as well the emergent algebra sector, with it's "emergent algebra crossing macro" described in section \ref{emergentalgebramoves}. (In fact this double appearance of knot diagram crossings is the path to the discovery of the extended graphic beta move.)

\paragraph{The moves.}         

\begin{enumerate}
\item[-]  \href{http://chorasimilarity.wordpress.com/2012/12/13/the-graphic-beta-move-with-details/}{The graphic beta move, with details}, section \ref{graphicbetawithdetails},   
\item[-]       \href{http://chorasimilarity.wordpress.com/2012/12/17/fan-out-moves-co-comm-co-assoc-global-fan-out-local-fan-out/}{Fan-out moves: CO-COMM, CO-ASSOC, GLOBAL FAN-OUT, LOCAL FAN-OUT}, section 
\ref{fanoutmoves},   
\item[-]  \href{http://chorasimilarity.wordpress.com/2012/12/17/pruning-moves/}{Pruning moves} and elimination of loops, section \ref{pruningmoves},     
\item[-] \href{http://chorasimilarity.wordpress.com/2012/12/20/emergent-algebra-moves-r1a-r1b-and-ext2}{Emergent algebra moves R1a, R1b, R2 and ext2}, section \ref{emergentalgebramoves}
\item[-]       \href{http://chorasimilarity.wordpress.com/2013/01/02/extensionality-in-graphic-lambda-calculus/}{The ext1 move}, section \ref{extensionalityingraphiclambda} for extensionality.   
\item[-]       \href{http://chorasimilarity.wordpress.com/tag/extended-beta-move/}{The extended beta move}, introduced     \href{http://chorasimilarity.wordpress.com/2012/12/09/graphic-beta-move-extended-to-explore/}{here}, section \ref{graphicbetamoveextendedtoexplore}.   
\end{enumerate}

The generating moves, discussed     \href{http://chorasimilarity.wordpress.com/2012/11/21/ancient-turing-machines-i-the-three-moirai/}{here},   \href{http://chorasimilarity.wordpress.com/2012/11/23/the-three-moirai-continued/}{here} and     \href{http://chorasimilarity.wordpress.com/2012/12/10/lachesis-computation-and-desire-to-explore-ancient-turing-machines-iii/}{here} are under consideration for futher inclusion in the list of moves.   

\section{The set GRAPH}

This is the set of graphs which are subjected to the moves.  Any finite assembly of the following elementary graphs, called "gates" is a graph in $ GRAPH$.
\begin{enumerate}
\item[-]   The $ \lambda$ gate, which corresponds to the lambda abstraction operation from lambda calculus, see     \href{http://en.wikipedia.org/wiki/Lambda_calculus#Lambda_terms}{lambda terms}. It is this:   

\vspace{.5cm}   \centerline{\includegraphics[width=20mm]{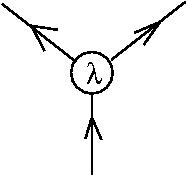}}  \vspace{.5cm}
\end{enumerate}
  
%    \href{http://chorasimilarity.wordpress.com/2012/12/12/introduction-to-graphic-lambda-calculus/lambdal-2/" rel="attachment wp-att-1502   <img class="aligncenter size-full wp-image-1502" alt="lambdal" src="http://chorasimilarity.files.wordpress.com/2012/12/lambdal1.jpg" width="105" height="100" />   

This gate looks like it has one input (the entry arrow) and two outputs (the left and right exit arrows respectively). This could not be a graph representing an operation, because an operation has two inputs and one output. For example, the lambda abstraction operation takes as inputs $ x$ a variable name and $ A$ a term and outputs the term $ \lambda x.A$.

Remember that the graphic lambda calculus does not have variable names.   There is a certain algorithm, described  \href{http://chorasimilarity.wordpress.com/2012/12/21/conversion-of-lambda-calculus-terms-into-graphs/}{here}, \href{http://arxiv.org/abs/1207.0332}{arXiv:1207.0332}, \cite{graphiclambdamoves} section 3 which  transforms a lambda term into a graph    in $ GRAPH$, such that to any lambda abstraction which appears in the term corresponds a $ \lambda$ gate. The algorithm starts with the representation of the lambda abstraction operation as a node with two inputs and one output, namely as an elementary gate which looks like the $ \lambda$  gate, but the orientation of the left exit arrow is inverse than the one of the $ \lambda$ gate. At some point in the algorithm the orientation is reversed and we get $ \lambda$ gates as shown here. 

The  $ \lambda$ gate looks like it takes a term $ A$ as input and it outputs at the left exit arrow the variable name $ x$ and at the right exit arrow the term $ \lambda x.A$. (It does not do this, properly, because there will be no variable names in the formalism.)
   
\begin{enumerate}
\item[-]  The $ \curlywedge$ graph, which corresponds to the application operation from lambda calculus, see     \href{http://en.wikipedia.org/wiki/Lambda_calculus#Lambda_terms}{lambda terms}. It is this:   

\vspace{.5cm}   \centerline{\includegraphics[width=20mm]{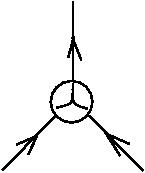}}  \vspace{.5cm}
\end{enumerate}
   
%    \href{http://chorasimilarity.wordpress.com/2012/12/12/introduction-to-graphic-lambda-calculus/curlyl/" rel="attachment wp-att-1507   <img class="aligncenter size-full wp-image-1507" alt="curlyl" src="http://chorasimilarity.files.wordpress.com/2012/12/curlyl.jpg" width="103" height="123" />   

This gate looks like the graph of an operation. The sign I use is like a curly join sign.
   
\begin{enumerate}
\item[-]   The $ \Upsilon$ graph, which will be used as a FAN-OUT gate, it is:   
 
  \vspace{.5cm}   \centerline{\includegraphics[width=20mm]{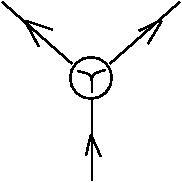}}  \vspace{.5cm}
\end{enumerate}

%    \href{http://chorasimilarity.wordpress.com/2012/12/12/introduction-to-graphic-lambda-calculus/upsil/" rel="attachment wp-att-1508   <img class="aligncenter size-full wp-image-1508" alt="upsil" src="http://chorasimilarity.files.wordpress.com/2012/12/upsil.jpg" width="102" height="103" />   
   
\begin{enumerate}
\item[-]     The $ \bar{\varepsilon}$ graph.       For any element $ \varepsilon \in \Gamma$ of an abelian group $ \Gamma$ (think about $ \Gamma$ as being $ (\mathbb{Z}, +)$ or $ ((0,+\infty), \cdot)$ ) there is an "exploration gate", or "dilation gate", which looks like the graph of an operation:   

 \vspace{.5cm}   \centerline{\includegraphics[width=20mm]{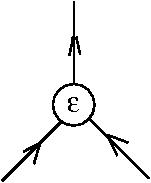}}  \vspace{.5cm}
\end{enumerate}
   
%    \href{http://chorasimilarity.wordpress.com/2012/12/12/introduction-to-graphic-lambda-calculus/epsil/" rel="attachment wp-att-1509   <img class="aligncenter size-full wp-image-1509" alt="epsil" src="http://chorasimilarity.files.wordpress.com/2012/12/epsil.jpg" width="102" height="124" />   

(Therefore we have a family of operations, called "dilations", indexed by the elements of an abelian group. This is a structure coming from  emergent algebras  \href{http://arxiv.org/abs/0907.1520}{arXiv:0907.1520}, \cite{buligairq}.)

We use these elementary graphs for constructing the graphs in $ GRAPH$. Any assembly of these gates, in any  (finite) number, which respects the orientation of arrows, is in $ GRAPH$.

Remark that we obtain trivalent graphs, with decorated nodes, each node having a cyclical order of his arrows (hence locally planar graphs).

Moreover, we accept arrows which input or output into nothing. Indeed, in particular the elementary graphs or gates are in $ GRAPH$ and all the arrows of an elementary graph either input or output to nothing.

Technically, we may imagine that we complete a graph in $ GRAPH$, if necessary, with univalent nodes, called "leaves" (they may be be decorated with "INPUT" or "OUTPUT", depending on the orientation of the arrow where they sit onto).
   
 For this reason we admit into $ GRAPH$ arrows without nodes which are elementary graphs, called       wires,             \includegraphics[width=20mm]{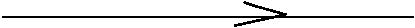} , and       loops       (without nodes from the elementary graphs, nor leaves)

 \vspace{.5cm}   \centerline{\includegraphics[width=20mm]{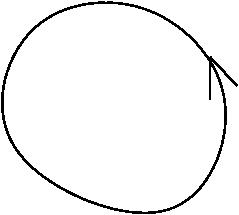}}  \vspace{.5cm}

%    \href{http://chorasimilarity.wordpress.com/2012/12/12/introduction-to-graphic-lambda-calculus/linel/" rel="attachment wp-att-1510   <img class="aligncenter size-full wp-image-1510" alt="linel" src="http://chorasimilarity.files.wordpress.com/2012/12/linel.jpg" width="267" height="18" />   

%    \href{http://chorasimilarity.wordpress.com/2012/12/12/introduction-to-graphic-lambda-calculus/loopl/" rel="attachment wp-att-1511   <img class="aligncenter size-full wp-image-1511" alt="loopl" src="http://chorasimilarity.files.wordpress.com/2012/12/loopl.jpg" width="144" height="129" />   
   
   Finally, we introduce an univalent gate,       the termination gate:    \includegraphics[width=5mm]{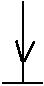} . 

%    \href{http://chorasimilarity.wordpress.com/2012/12/12/introduction-to-graphic-lambda-calculus/topl/" rel="attachment wp-att-1512   <img class="aligncenter size-full wp-image-1512" alt="topl" src="http://chorasimilarity.files.wordpress.com/2012/12/topl.jpg" width="44" height="86" />   

The termination gate has an input leaf and no output.

Any graph which is a reunion of lines, loops and assemblies of the elementary graphs (termination graph included) is in $ GRAPH$ and any graph in $GRAPH$ is obtained in this way.

\section{The graphic beta move, with details}
\label{graphicbetawithdetails}

Here it is explained the most important move in graphic lambda calculus: the graphic beta move.

A move is a transformation of a graph into another. It is not a function from $ GRAPH$ to $ GRAPH$ though, as we shall see.

The graphic beta move takes its name from the     \href{http://en.wikipedia.org/wiki/Lambda_calculus#Beta_reduction}{beta reduction in lambda calculus}.

Consider the set $ GRAPH$ with the equality of oriented,  locally planar graphs (with leaves not numbered), denoted by $ \equiv$.

A graphic beta move is a local move which consists in replacing a pattern (sub-graph)  by another pattern. More precisely, we can apply the graphic beta move to any graph which contains a certain pattern (sub-graph) and we obtain a new graph, which is identical with the initial one, with the exception of the pattern, now replaced by the other pattern.  In some figures, for clarity, we shall add a dashed closed curve which encircles the pattern which is subjected to the move. This is the convention followed in 
 \href{http://arxiv.org/abs/1207.0332}{arXiv:1207.0332}, \cite{graphiclambdamoves}. 

\begin{enumerate}
\item[-] The graphic beta move. See the explanations after the figure.
\end{enumerate}

\vspace{.5cm}   \centerline{ \includegraphics[width=100mm]{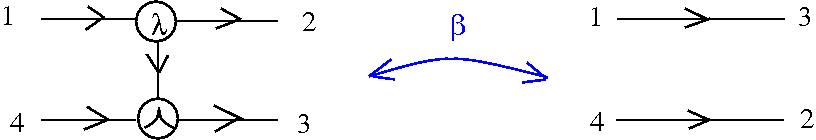}}  \vspace{.5cm}

%    \href{http://chorasimilarity.wordpress.com/2012/12/13/the-graphic-beta-move-with-details/beta_move_1/" rel="attachment wp-att-1515   <img class="aligncenter size-full wp-image-1515" alt="beta_move_1" src="http://chorasimilarity.files.wordpress.com/2012/12/beta_move_1.jpg" width="595" height="102" />   

The labels "1, 2, 3, 4" are used only as guides for gluing correctly the new pattern, after removing the old one. Imagine that we have a graph which contains as a sub-graph the one from the left hand side (LHS) of the previous figure. Then, by a graphic beta move, we may replace the pair of nodes from the LHS figure with two arrows, as indicated by the labelling "1,2,3,4". The rest of the graph is unchanged.

We may represent the graphic beta move like this:

\vspace{.5cm}   \centerline{ \includegraphics[width=100mm]{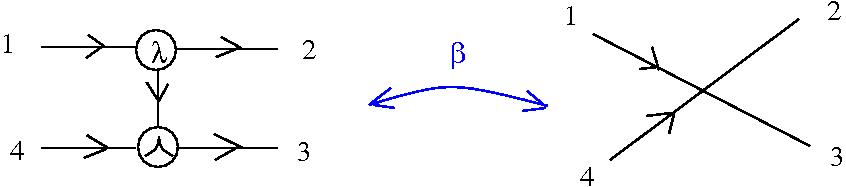}}  \vspace{.5cm}

%    \href{http://chorasimilarity.wordpress.com/2012/12/13/the-graphic-beta-move-with-details/beta_move_6/" rel="attachment wp-att-1528   <img class="aligncenter size-full wp-image-1528" alt="beta_move_6" src="http://chorasimilarity.files.wordpress.com/2012/12/beta_move_6.jpg" width="595" height="131" />   

Up to the particular embedding in the plane of the graph from the right hand side (RHS), the move is the same. The intersection of the "1,3" arrow with the "4,2" arrow is an artifact of the embedding, there is no node there. Intersections of arrows have no meaning, remember that we work with graphs which are locally planar, not globally planar.

The graphic beta move goes into both directions. Referring to the first two pictures,  we may pick a pair of arrows and label them with "1,2,3,4", such that, according to the orientation of the arrows,  "1" points to "3" and "4" points to "2", without any node or label between "1" and "3" and between "4" and "2" respectively. Then, by a graphic beta move, we may replace the portions of the two arrows which are between "1" and "3", respectively between "4" and "2", by the pattern from the LHS of the figure.

Here is an example.

 \vspace{.5cm}   \centerline{\includegraphics[width=135mm]{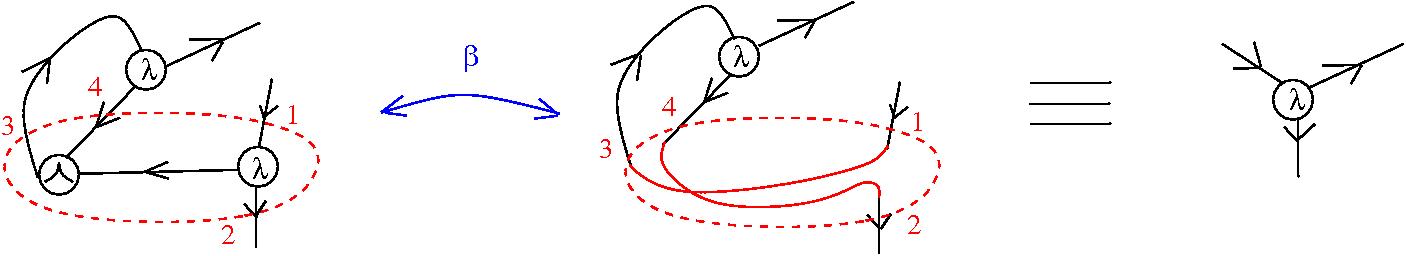}}  \vspace{.5cm}

%    \href{http://chorasimilarity.wordpress.com/2012/12/13/the-graphic-beta-move-with-details/beta_move_2/" rel="attachment wp-att-1516   <img class="aligncenter size-full wp-image-1516" alt="beta_move_2" src="http://chorasimilarity.files.wordpress.com/2012/12/beta_move_2.jpg" width="637" height="115" />   

The red drawings have no meaning in the formalism, I have put them only to make a more clear representation  of the move. Remark that we may read the move from left to right, as well as from right to left. The dashed closed curve serves to mark the pattern.

The graphic beta move may be applied to a single arrow, or to a loop, there is nothing which might stop this. For example, here are three applications of the graphic beta move:

\vspace{.5cm}   \centerline{ \includegraphics[width=120mm]{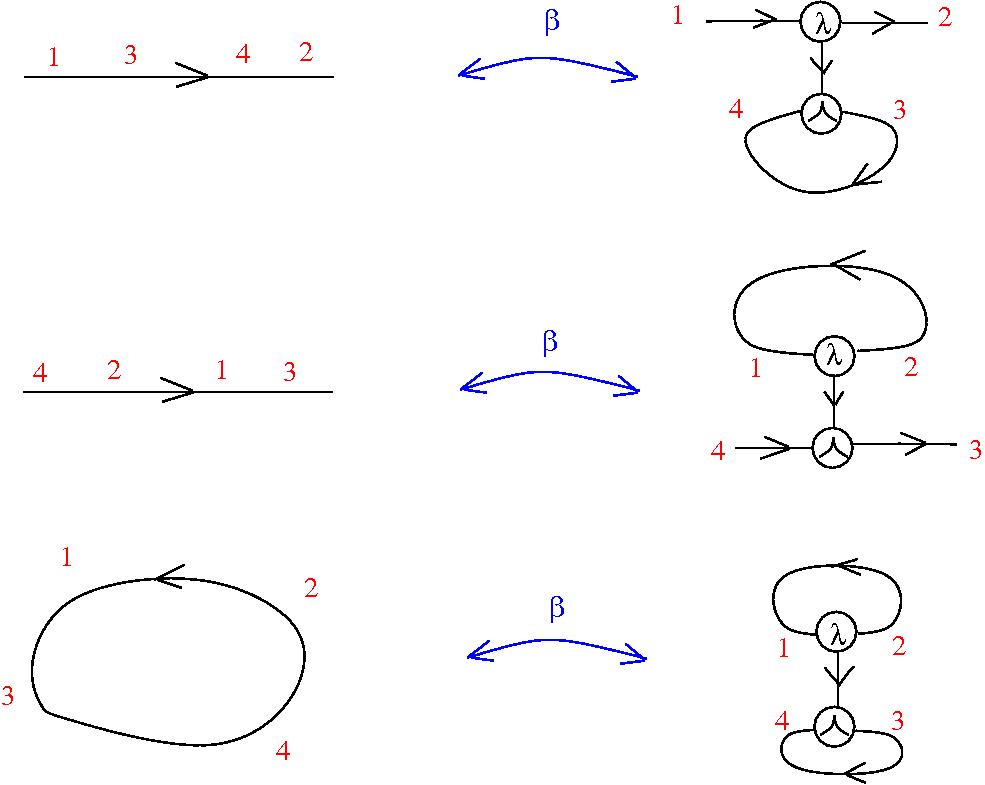}}  \vspace{.5cm}

%    \href{http://chorasimilarity.wordpress.com/2012/12/13/the-graphic-beta-move-with-details/beta_move_3/" rel="attachment wp-att-1517   <img class="aligncenter size-full wp-image-1517" alt="beta_move_3" src="http://chorasimilarity.files.wordpress.com/2012/12/beta_move_3.jpg" width="595" height="474" />   

In particular this gives a reason for considering loops and wires as members of $ GRAPH$.

A point important to stress is that the graphic beta move is not dependent of the particular embedding into the plane of the graphs (which we use when we draw the graphs). This will be important further, when we shall  discuss about the crossing macros, see  \href{http://chorasimilarity.wordpress.com/2012/09/06/graphic-beta-rule-as-braiding/}{Graphic beta rule as braiding} section \ref{graphicbetaruleasbraiding}.

Also, different labelling (of the same graphs) leads to different outcomes of the graphic beta move. This is explained in the next figure: 

\vspace{.5cm}   \centerline{ \includegraphics[width=135mm]{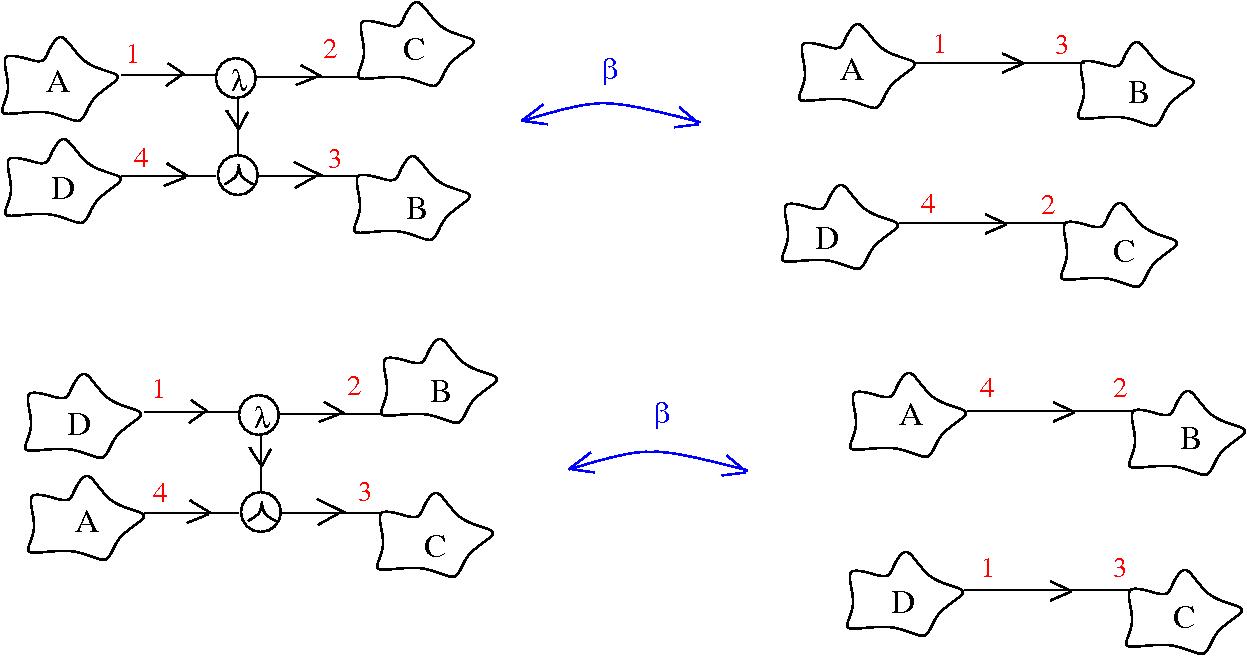}}  \vspace{.5cm}

%    \href{http://chorasimilarity.wordpress.com/2012/12/13/the-graphic-beta-move-with-details/beta_move_4-2/" rel="attachment wp-att-1523   <img class="aligncenter size-full wp-image-1523" alt="beta_move_4" src="http://chorasimilarity.files.wordpress.com/2012/12/beta_move_41.jpg" width="595" height="313" />   

Here "A, B, C, D" denote four graphs in $ GRAPH$. Two possible graphic beta moves are depicted. Let's read them from right to left.

At the right we have a graph with two arrows, connecting "A" with "B" and "D" with "C". We may apply the graphic beta move to this graph in two ways, depending on the labelling of the arrows. If we mark the arrow from "A" to "B" with "1,3" and the arrow from "D" to "C" with "4,2" then we get the first graphic beta move. However, we may choose to mark the arrow from "A" to "B" with "4,2" and the arrow from "D" to "C" with "1,3", applying then the graphic beta move. The outcome is different.

A particular case of the previous figure is yet another justification for having loops as elements in $ GRAPH$.

\vspace{.5cm}   \centerline{\includegraphics[width=135mm]{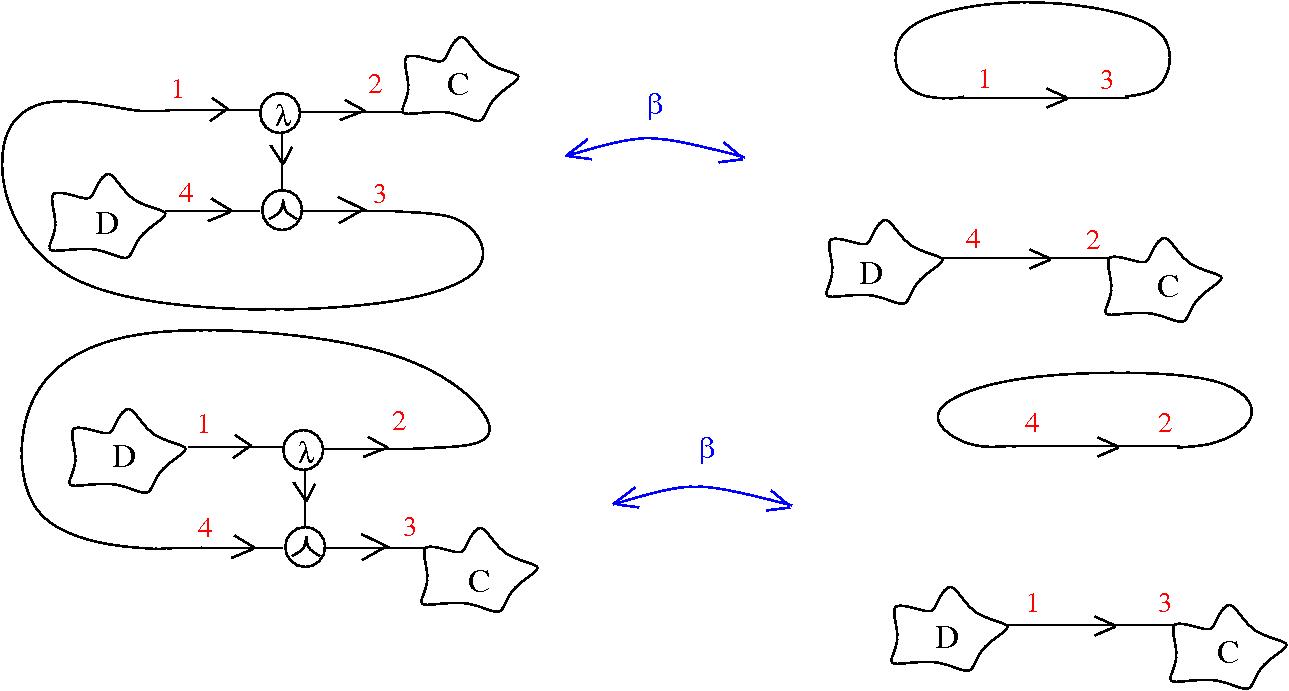}}  \vspace{.5cm}

%    \href{http://chorasimilarity.wordpress.com/2012/12/13/the-graphic-beta-move-with-details/beta_move_5/" rel="attachment wp-att-1525   <img class="aligncenter size-full wp-image-1525" alt="beta_move_5" src="http://chorasimilarity.files.wordpress.com/2012/12/beta_move_5.jpg" width="595" height="318" />   

These two applications of the graphic beta move may be represented alternatively like this:

\vspace{.5cm}   \centerline{\includegraphics[width=135mm]{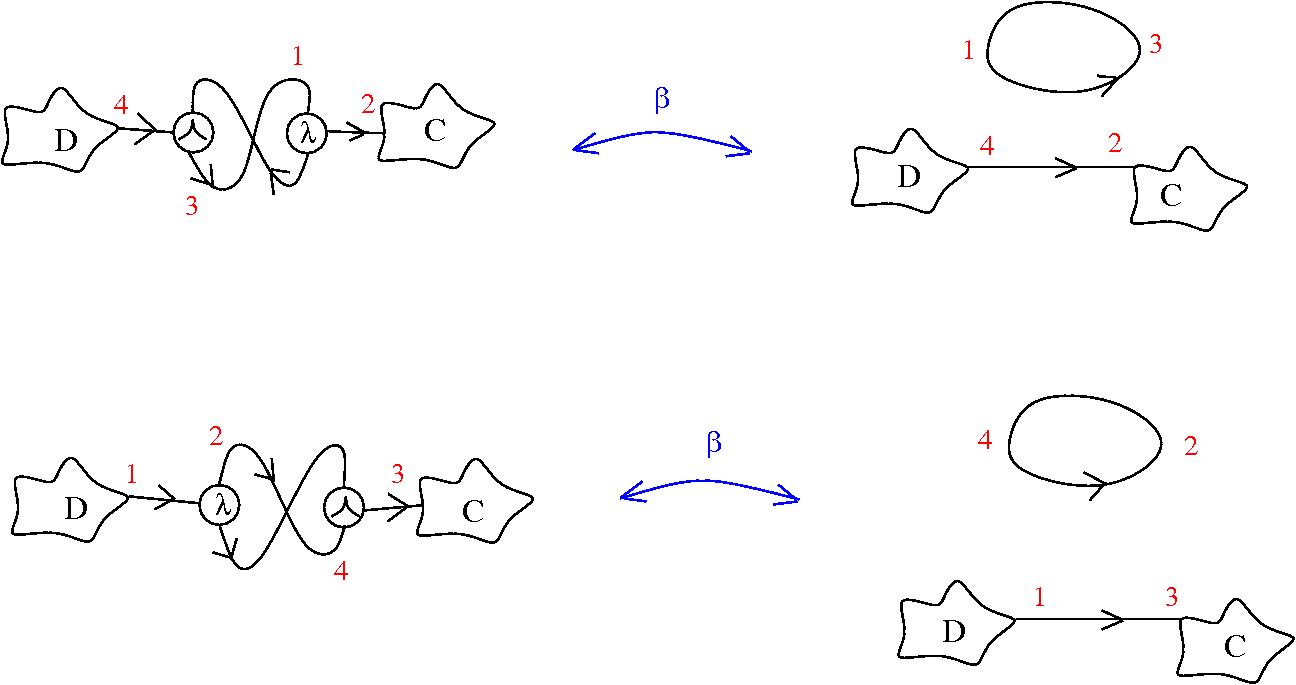}}  \vspace{.5cm}

%    \href{http://chorasimilarity.wordpress.com/2012/12/13/the-graphic-beta-move-with-details/beta_move_8/" rel="attachment wp-att-1531   <img class="aligncenter size-full wp-image-1531" alt="beta_move_8" src="http://chorasimilarity.files.wordpress.com/2012/12/beta_move_8.jpg" width="595" height="314" />   

A move which looks very much alike   the graphic beta move is the   UNZIP operation   from  the formalism of  knotted trivalent graphs, see for example the paper      \href{http://arxiv.org/abs/math/0311458}{The algebra of knotted trivalent graphs and Turaev’s shadow world}    by D.P. Thurston, section 3. 
In order to see this, let's draw again the graphic beta move, this time without labelling the arrows:

\vspace{.5cm}   \centerline{\includegraphics[width=120mm]{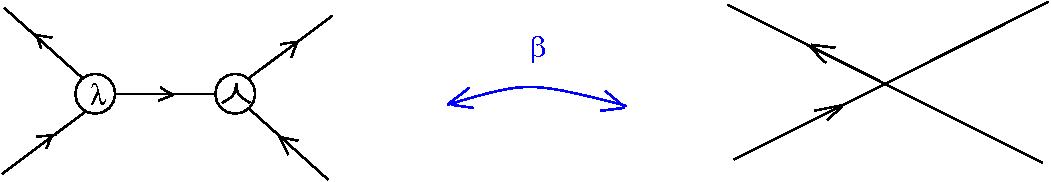}}  \vspace{.5cm}

%    \href{http://chorasimilarity.wordpress.com/2012/12/13/the-graphic-beta-move-with-details/beta_move_9-2/" rel="attachment wp-att-1534   <img class="aligncenter size-full wp-image-1534" alt="beta_move_9" src="http://chorasimilarity.files.wordpress.com/2012/12/beta_move_91.jpg" width="595" height="103" />   

The unzip operation acts only from left to right in the following figure. Remarkably, it acts on trivalent graphs (but not oriented).

\vspace{.5cm}   \centerline{\includegraphics[width=120mm]{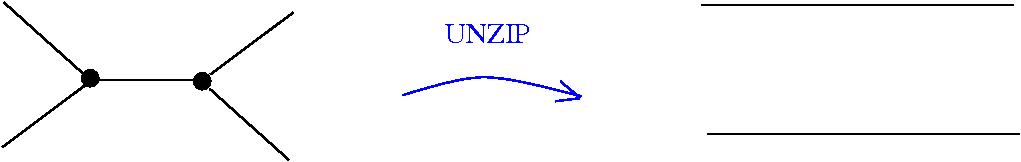}}  \vspace{.5cm}

%    \href{http://chorasimilarity.wordpress.com/2012/12/13/the-graphic-beta-move-with-details/beta_move_10/" rel="attachment wp-att-1535   <img class="aligncenter size-full wp-image-1535" alt="beta_move_10" src="http://chorasimilarity.files.wordpress.com/2012/12/beta_move_10.jpg" width="595" height="94" />   

UNZIP is not a move in graphic lambda calculus.

\section{Fan-out moves: CO-COMM, CO-ASSOC, GLOBAL FAN-OUT, LOCAL FAN-OUT}
\label{fanoutmoves}

Here are described the moves directly related to the $ \Upsilon$ gate.

\begin{enumerate}
\item[-]     CO-COMM move.              This is the local move depicted in the following figure   

\vspace{.5cm}   \centerline{\includegraphics[width=80mm]{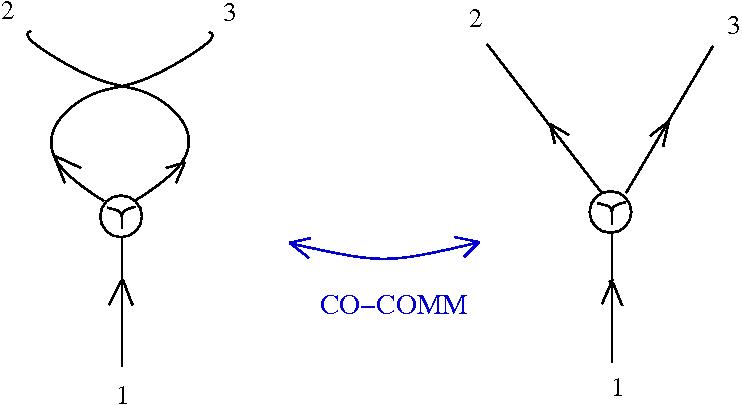}}  \vspace{.5cm}  
\end{enumerate}

%    \href{http://chorasimilarity.wordpress.com/2012/12/17/fan-out-moves-co-comm-co-assoc-global-fan-out-local-fan-out/commr/" rel="attachment wp-att-1558   <img class="aligncenter size-full wp-image-1558" alt="commr" src="http://chorasimilarity.files.wordpress.com/2012/12/commr.jpg" width="595" height="324" />   

It means we may permute the outputs of a $ \Upsilon$ gate. The name means "co-commutativity", because the diagram resembles to the one of the commutativity property, with the exception of the arrows orientations, which are in opposite directions (hence "co-").
   
\begin{enumerate}
\item[-]         CO-ASSOC move.       This is a local move, described by the next figure   

\vspace{.5cm}   \centerline{\includegraphics[width=80mm]{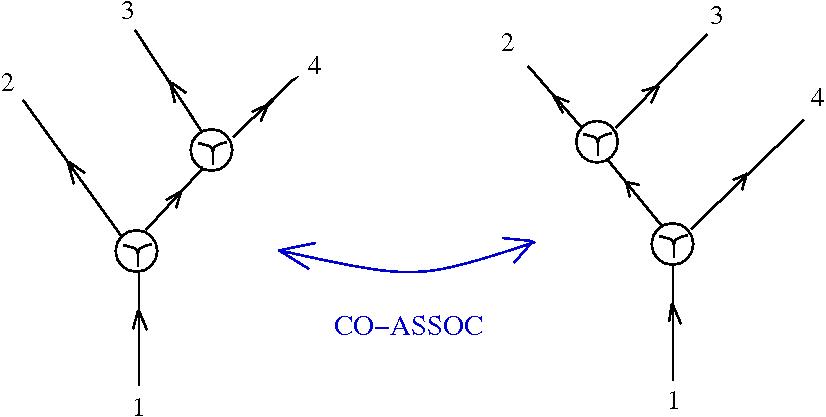}}  \vspace{.5cm}  

\end{enumerate}

%    \href{http://chorasimilarity.wordpress.com/2012/12/17/fan-out-moves-co-comm-co-assoc-global-fan-out-local-fan-out/assocr-2/" rel="attachment wp-att-1559   <img class="aligncenter size-full wp-image-1559" alt="assocr" src="http://chorasimilarity.files.wordpress.com/2012/12/assocr.jpg" width="595" height="300" />   

It has the following effect: by using CO-ASSOC moves, we may move between any two binary trees formed only with $ \Upsilon$ gates, with the same number of output leaves. The name means "co-associativity" and the explanation is similar to the previous one, with "commutativity" replaced by "associativity".
 
\begin{enumerate}
\item[-]  GLOBAL FAN-OUT.       This is a global move, because it involves a modification of an arbitrary number of nodes (gates) and arrows.   
\end{enumerate}
   
Precisely, the move acts like in the following picture: if $ A$ is a graph in $ GRAPH$  then we may replace the graph $ A$ connected to an $ \Upsilon$  gate by two copies of $ A$. More specifically, the graph $A$ from the left hand side of the picture has to be connected only to the input of the $\Upsilon$ gate, with no other connection to the rest of the graph which may be found outside the figured dashed curve. 

\vspace{.5cm}   \centerline{\includegraphics[width=100mm]{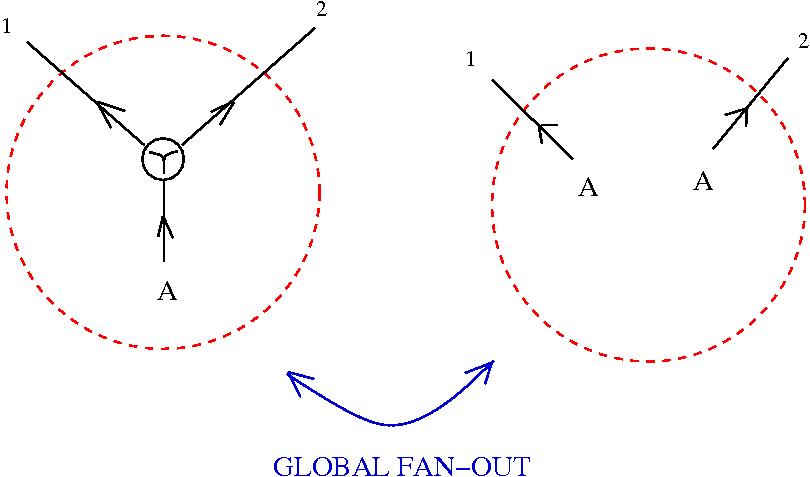}}  \vspace{.5cm}

%    \href{http://chorasimilarity.wordpress.com/2012/12/17/fan-out-moves-co-comm-co-assoc-global-fan-out-local-fan-out/fanoutr/" rel="attachment wp-att-1560   <img class="aligncenter size-full wp-image-1560" alt="fanoutr" src="http://chorasimilarity.files.wordpress.com/2012/12/fanoutr.jpg" width="595" height="350" />   

GLOBAL FAN-OUT implies CO-COMM  (namely two GLOBAL FAN-OUT moves have the effect of one CO-COMM move).  The move CO-COMM is not useless though: we may choose to work in a sector of the graphic lambda calculus which uses CO-COMM but not GLOBAL FAN-OUT.

There is a variant of this move which is local.
\begin{enumerate}   
\item[-]         LOCAL FAN-OUT.       Fix a number $ N$ and consider only graphs $ A$ which have at most $ N$ (nodes + arrows). The $ N$ LOCAL FAN-OUT move is the same as the GLOBAL FAN-OUT move, only it applies only to such graphs $ A$.   
\end{enumerate}
   
LOCAL FAN-OUT does not imply CO-COMM.

\paragraph{Important remark:}    The gate $ \Upsilon$ is really a fan-out gate only in the sense described by the FAN-OUT moves. In the absence of one of this moves, the gate cannot  be described as a "fan-out". The "fan-out" is in the moves, not in the gate.

\section{Pruning moves. Elimination of loops}
\label{pruningmoves}

Here are described the moves related to the univalent termination gate. There are local and global pruning moves. At the end of the section we write about the elimination of loops. 

\begin{enumerate}
\item[-]         LOCAL PRUNING.        The local moves which eliminate "dead edges" or "dead nodes" are the following:   
   
\vspace{.5cm}   \centerline{\includegraphics[width=80mm]{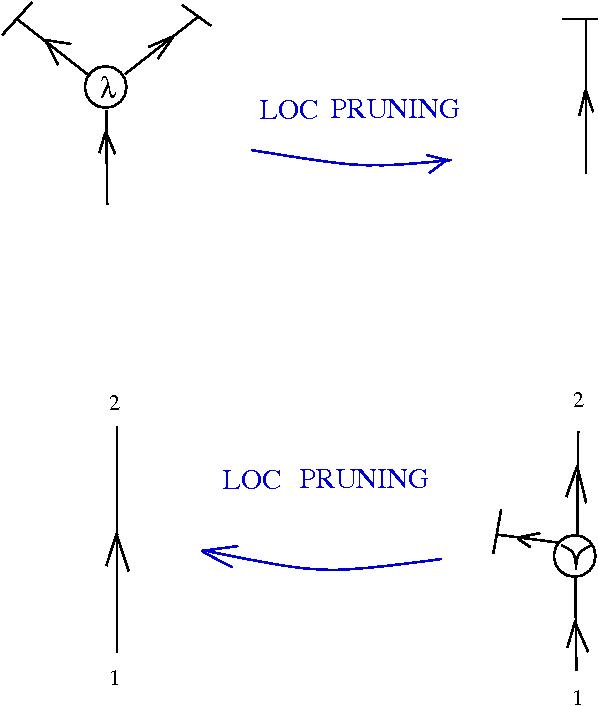}}  \vspace{.5cm}   

%    \href{http://chorasimilarity.wordpress.com/2012/12/17/pruning-moves/lambdapr-2/" rel="attachment wp-att-1579   <img class="aligncenter size-full wp-image-1579" alt="lambdapr" src="http://chorasimilarity.files.wordpress.com/2012/12/lambdapr1.jpg" width="425" height="499" />   

\vspace{.5cm}   \centerline{\includegraphics[width=100mm]{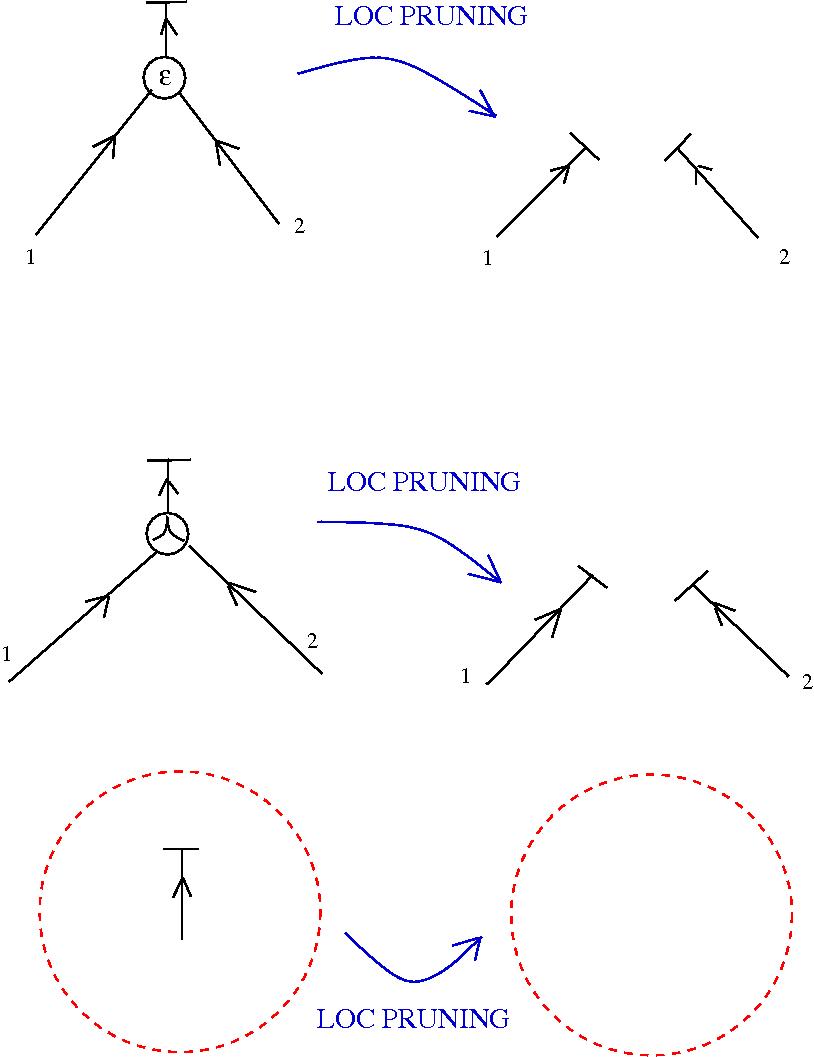}}  \vspace{.5cm}  

%    \href{http://chorasimilarity.wordpress.com/2012/12/17/pruning-moves/epsipr/" rel="attachment wp-att-1570   <img class="aligncenter size-full wp-image-1570" alt="epsipr" src="http://chorasimilarity.files.wordpress.com/2012/12/epsipr.jpg" width="531" height="688" />   

\end{enumerate}

These moves go in one direction, compared with the graphic beta move, which is allowed in both directions. In      \href{http://arxiv.org/abs/1207.0332v1}{arXiv:1207.0332v1}   \cite{graphiclambdamoves} paragraph 2.6, where this moves were introduced, they are allowed in both directions. However, I intend to modify this and to introduce generation moves instead, like CREA and GARB, as explained in Ancient Turing machines posts     \href{http://chorasimilarity.wordpress.com/2012/11/21/ancient-turing-machines-i-the-three-moirai/}{here}    and     \href{http://chorasimilarity.wordpress.com/2012/11/23/the-three-moirai-continued/}{here}.
   
\begin{enumerate}
\item[-]         GLOBAL PRUNING.       Like  the     \href{http://chorasimilarity.wordpress.com/2012/12/17/fan-out-moves-co-comm-co-assoc-global-fan-out-local-fan-out/}{GLOBAL FAN-OUT move} section \ref{fanoutmoves}, this is a global move. A local version of it, in the spirit of the     \href{http://chorasimilarity.wordpress.com/2012/12/17/fan-out-moves-co-comm-co-assoc-global-fan-out-local-fan-out/}{LOCAL FAN-OUT move}, may be introduced.  For any graph $ A$ in $ GRAPH$, if $ A$  is connected ONLY to a termination gate, then we may remove it.   

\vspace{.5cm}   \centerline{\includegraphics[width=100mm]{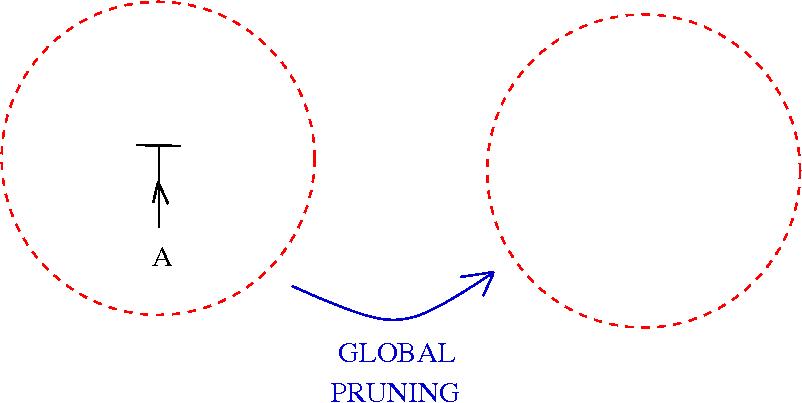}}  \vspace{.5cm}   
   
%    \href{http://chorasimilarity.wordpress.com/2012/12/17/pruning-moves/globalpr/" rel="attachment wp-att-1571   <img class="aligncenter size-full wp-image-1571" alt="globalpr" src="http://chorasimilarity.files.wordpress.com/2012/12/globalpr.jpg" width="595" height="298" />   

\end{enumerate}

\begin{enumerate}
\item[-] Elimination of loops. This move has been proposed in  \href{http://arxiv.org/abs/1207.0332}{arXiv:1207.0332}, \cite{graphiclambdamoves} section 2, paragraph 2.10. We are free to erase or add loops to graphs in $GRAPH$. 
\end{enumerate}

\section{Emergent algebra moves R1a, R1b, R2 and ext2}
\label{emergentalgebramoves}

Here are described the moves related to emergent algebras. This moves involve the gates $ \bar{\varepsilon}$, $ \Upsilon$ and the termination gate.

I shall use the notations of the Reidemeister moves for oriented knot diagrams from the paper  \href{http://arxiv.org/abs/0908.3127}{arXiv:0908.3127} \cite{polyak},   only that I use the letter “R” from “Reidemeister” instead of “$ \Omega$” used by Polyak. 

In the setting of emergent algebras \href{http://arxiv.org/abs/1005.5031}{arXiv:1005.5031} \cite{buligabraided}, \href{http://arxiv.org/abs/0907.1520}{arXiv:0907.1520}, \cite{buligairq}, to each $\varepsilon \in \Gamma$, where $\Gamma$ is an abelian group, is associated an operation (on a set $X$), which is graphically represented (in the graphic lambda calculus formalism) by an associated gate $ \bar{\varepsilon}$. With respect to any  $\varepsilon \in \Gamma$ the pair $(X, \bar{\varepsilon})$ is an idempotent right quasigroup, moreover the operation associated to the neutral element of $\Gamma$ is $\displaystyle x \bar{1} y = y$ and for any $\varepsilon, \mu \in \Gamma$ and any $x, y \in X$ we have $\displaystyle x \bar{\varepsilon} \left( x \bar{\mu} y \right) = x \, \bar{\varepsilon\mu} \, y$. There are special cases, namely when the operations $\displaystyle \bar{\varepsilon}$ are quasigroups operations, and the most particular case of left distributive emergent algebras, where in particular $\displaystyle (X, \bar{\varepsilon})$ becomes a quandle. 

These algebraic axioms admit a description in terms of knot diagrams (via the well-known relation between quandles and knot diagrams). In the paper \href{http://arxiv.org/abs/1005.5031}{arXiv:1005.5031} \cite{buligabraided}, then in  \href{http://arxiv.org/abs/1103.6007}{arXiv:1103.6007}, \cite{buligachora} section 3.1 this is explained in detail.  

In the decorated knots formalism, crossings of oriented wires are decorated with elements $ \varepsilon$ of a commutative group $ \Gamma$. The relation between these crossings and their representations in terms of trivalent graphs is described by the following figures, which can be seen as definitions of the (emergent algebra) "crossing macros" (see section \ref{graphicbetaruleasbraiding} for the definition of lambda calculus crossing macros): 

\vspace{.5cm}   \centerline{\includegraphics[width=100mm]{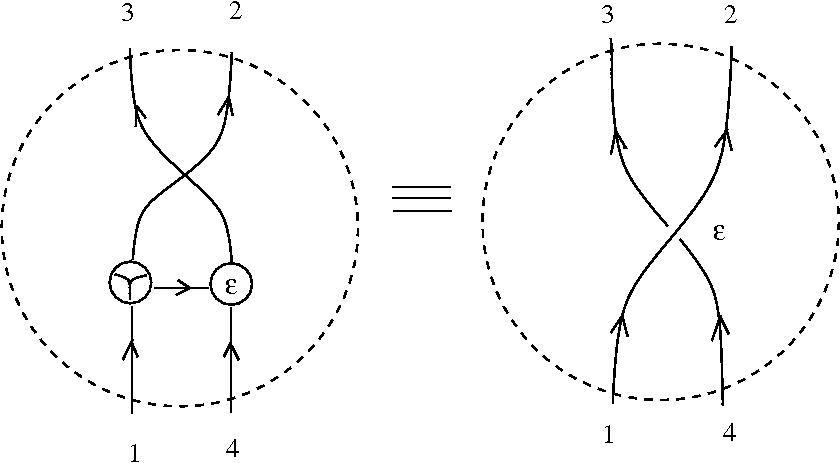}}  \vspace{.5cm}  

%    \href{http://chorasimilarity.files.wordpress.com/2012/10/emerr_dif_11.jpg   <img class="aligncenter  wp-image-1182" title="emerr_dif_11" alt="" src="http://chorasimilarity.files.wordpress.com/2012/10/emerr_dif_11.jpg?w=300" height="257" width="468" />   

\vspace{.5cm}   \centerline{\includegraphics[width=100mm]{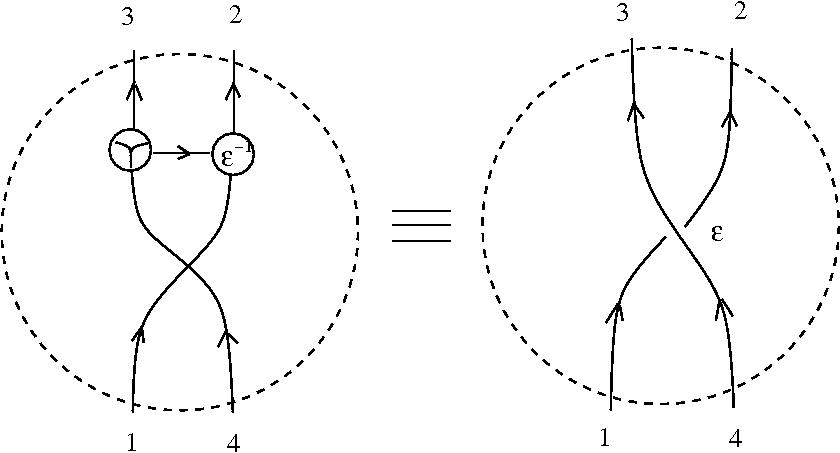}}  \vspace{.5cm}

%    \href{http://chorasimilarity.files.wordpress.com/2012/10/emerr_dif_21.jpg   <img class="aligncenter  wp-image-1183" title="emerr_dif_21" alt="" src="http://chorasimilarity.files.wordpress.com/2012/10/emerr_dif_21.jpg?w=300" height="259" width="483" />   

(See also the post  \href{http://chorasimilarity.wordpress.com/2012/10/31/3d-crossings-in-emergent-algebras/}{3D crossings in emergent algebras} in order to understand the relation with knot diagram crossings and with Reidemeister moves for oriented knot diagrams.)

The emergent algebra moves are just imported from the emergent algebra formalism. They have the following form. 

\begin{enumerate}
\item[-]         The move R1a       is described in the following figure   
   
\vspace{.5cm}   \centerline{\includegraphics[width=100mm]{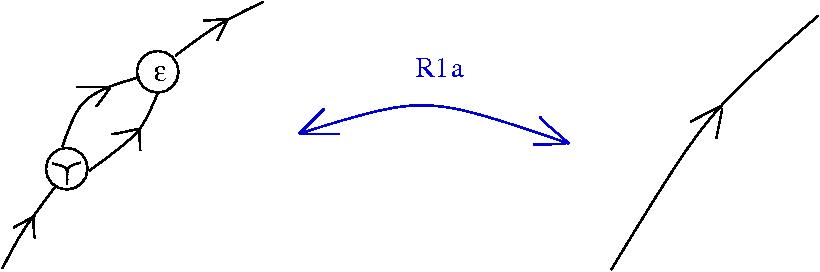}}  \vspace{.5cm} 

\end{enumerate}

%    \href{http://chorasimilarity.wordpress.com/2012/12/20/emergent-algebra-moves-r1a-r1b-and-ext2/r1amove/" rel="attachment wp-att-1620   <img class="aligncenter size-full wp-image-1620" alt="r1amove" src="http://chorasimilarity.files.wordpress.com/2012/12/r1amove.jpg" width="595" height="197" />   

The reason for calling this move R1a is that is related to the move R1a from Polyak paper ( and also to R1d move). In the papers on emergent algebras this move is called R1, that is "Reidemeister move 1". (In order to obtain exactly the Reidemeister move R1a from the oriented knot diagrams formalism we need also CO-COMM.)

\begin{enumerate}
\item[-]         The move R1b       is this:   
   
\vspace{.5cm}   \centerline{\includegraphics[width=100mm]{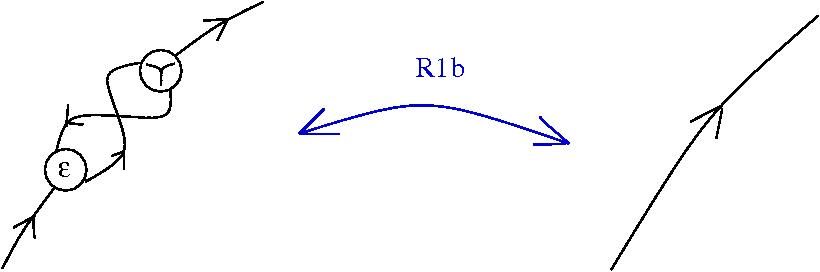}}  \vspace{.5cm}  

\end{enumerate}

%    \href{http://chorasimilarity.wordpress.com/2012/12/20/emergent-algebra-moves-r1a-r1b-and-ext2/r1bmove/" rel="attachment wp-att-1621   <img class="aligncenter size-full wp-image-1621" alt="r1bmove" src="http://chorasimilarity.files.wordpress.com/2012/12/r1bmove.jpg" width="595" height="197" />   

This move is related to the move R1b from Polyak (and also to R1c).   This move does not appear in relation with general emergent algebras, it is true only for a special subclass of them, namely uniform idempotent quasigroups.       
   
\begin{enumerate}
\item[-]         The move R2        is the following:   
   
\vspace{.5cm}   \centerline{\includegraphics[width=100mm]{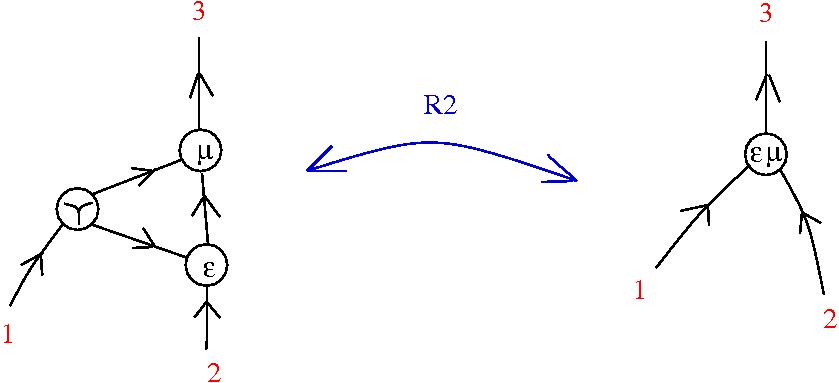}}  \vspace{.5cm}  

\end{enumerate}

%    \href{http://chorasimilarity.wordpress.com/2012/12/20/emergent-algebra-moves-r1a-r1b-and-ext2/r2move/" rel="attachment wp-att-1666   <img class="aligncenter size-full wp-image-1666" alt="r2move" src="http://chorasimilarity.files.wordpress.com/2012/12/r2move1.jpg" width="595" height="271" />   

This move is related to all Reidemeister 2 moves (using also CO-ASSOC, CO-COMM, and LOCAL PRUNING).

\begin{enumerate}
\item[-]         The move ext2.       The notation comes from the rule (ext2) from lambda-Scale calculus. In emergent algebra language it means that the emergent algebra operation indexed by the neutral element of $1 \in \Gamma$ is the trivial operation $ x \bar{1} y = y$.   
   
\vspace{.5cm}   \centerline{\includegraphics[width=100mm]{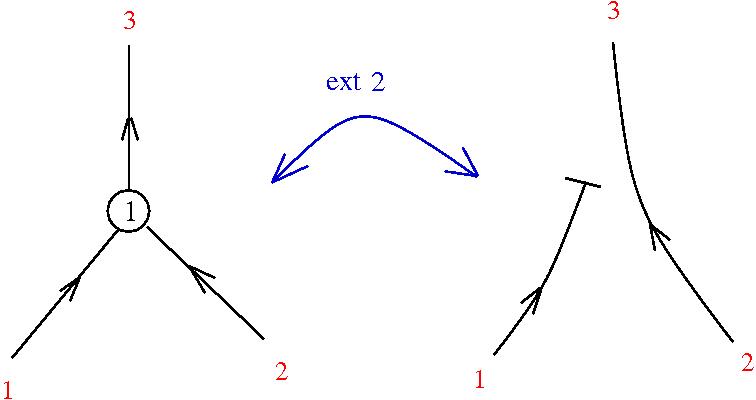}}  \vspace{.5cm} 

\end{enumerate}

%    \href{http://chorasimilarity.wordpress.com/2012/12/20/emergent-algebra-moves-r1a-r1b-and-ext2/ext2r/" rel="attachment wp-att-1623   <img class="aligncenter size-full wp-image-1623" alt="ext2r" src="http://chorasimilarity.files.wordpress.com/2012/12/ext2r.jpg" width="595" height="314" />   

(but for this LOCAL PRUNING is also used).

\section{Extensionality in graphic lambda calculus}
\label{extensionalityingraphiclambda}

I want to discuss here the introduction of extensionality in the graphic lambda calculus. In some sense, extensionality is already present in the \href{http://chorasimilarity.wordpress.com/2012/12/20/emergent-algebra-moves-r1a-r1b-and-ext2/}{emergent algebra moves} section \ref{emergentalgebramoves}. Have you noticed the (name of the) move "ext2"?

However, the \href{http://en.wikipedia.org/wiki/Lambda_calculus#.CE.B7-conversion}{eta-reduction}      from untyped lambda calculus needs a new move. I called it the ext1 move in \href{http://arxiv.org/abs/1207.0332}{arXiv:1207.0332}  \cite{graphiclambdamoves} paragraph 2.7. It is a global move, because in order to use it one has to check a global condition (without bound on the number of nodes and edges involved in the condition). In the mentioned paper I stated that the move applies only to graphs  \href{http://chorasimilarity.wordpress.com/2012/12/21/conversion-of-lambda-calculus-terms-into-graphs/}{which represent lambda calculus terms}, but now I see no reason why it has to be confined to this sector, therefore here I shall formulate the ext1 move in more generality.

\begin{enumerate}
\item[-] The move ext1.    If there is no oriented path from "2" to "1" outside the left hand side picture then one may replace this picture by an edge. Conversely, if there is no oriented path connecting "2" with "1" then one may replace the edge with the graph from the left hand side of the following picture:
\end{enumerate}

\vspace{.5cm}   \centerline{\includegraphics[width=80mm]{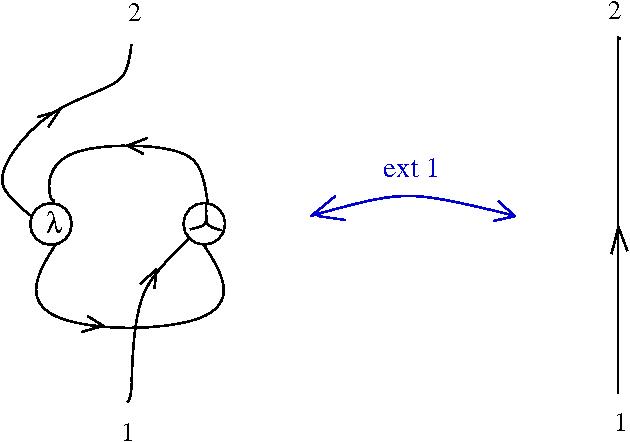}}  \vspace{.5cm} 

%\href{http://chorasimilarity.wordpress.com/2013/01/02/extensionality-in-graphic-lambda-calculus/ext1r/" rel="attachment wp-att-1713"><img class="aligncenter size-full wp-image-1713" alt="ext1r" src="http://chorasimilarity.files.wordpress.com/2013/01/ext1r1.jpg" width="595" height="417" />     

This move acts like eta-reduction when translated back from graphs to lambda calculus terms. "Ext" comes from "extensionality".

Let us see why we need to formulate the move like this. Suppose that we eliminate the global condition  "there is no oriented path from "2" to "1" outside the left hand side  of the previous picture". In particular let us suppose that there is an edge from "2" to "1" which completes the graphs from the LHS of the previous picture. For this graph, which appears at left in the next figure, we may use the graphic beta move like this (numbers in red indicate how we use the graphic beta move):

\vspace{.5cm}   \centerline{\includegraphics[width=80mm]{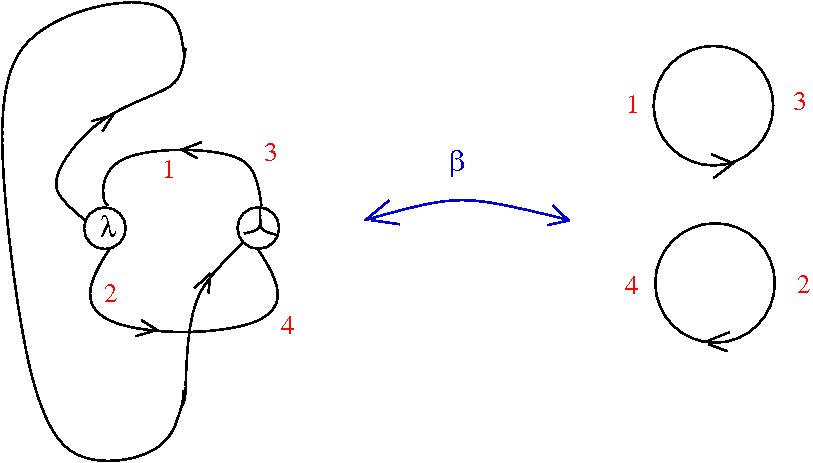}}  \vspace{.5cm} 

%\href{http://chorasimilarity.wordpress.com/2013/01/02/extensionality-in-graphic-lambda-calculus/non_ext1_move/" rel="attachment wp-att-1715"><img class="aligncenter size-full wp-image-1715" alt="non_ext1_move" src="http://chorasimilarity.files.wordpress.com/2013/01/non_ext1_move.jpg" width="595" height="338" />     

If we could apply the ext 1 move to this graph, then the result would be the following:

\vspace{.5cm}   \centerline{\includegraphics[width=80mm]{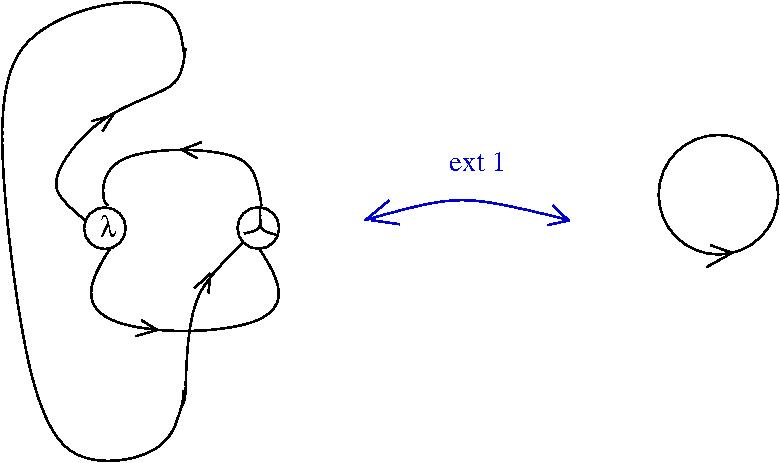}}  \vspace{.5cm}

%\href{http://chorasimilarity.wordpress.com/2013/01/02/extensionality-in-graphic-lambda-calculus/non_ext1_move2/" rel="attachment wp-att-1726"><img class="aligncenter size-full wp-image-1726" alt="non_ext1_move2" src="http://chorasimilarity.files.wordpress.com/2013/01/non_ext1_move21.jpg" width="595" height="353" />     

Not only the result of this move is one loop, but it is the "wrong" loop, in the sense that this loop is obtained by closing the edge which vanishes  when the graphic beta move is applied. There is no contradiction though, unless we wish to decorate the edges (i.e. to evaluate the result of the computation). It is just strange.

There is one question left: why the move called "ext2", which is one of the emergent algebra moves, is also an extensionality move? An answer is that the move ext2 says that the dilation of coefficient "1" is the identity  {\em  function}.

\section{Graphic beta rule as braiding}
\label{graphicbetaruleasbraiding}

The moves of the graphic lambda calculus are independent of the particular embeddings of trivalent graphs into the plane. This may not be obvious from the figures which explain the moves, but it is so after further inspection.

Let us look again at the graphic beta move. In fact, this is a move which transforms a pair of edges (from the right of the picture) into the graph from the left of the picture. The fact that the edges intersect in the figure is irrelevant. What matters is that, for the sake of performing the move, one edge is temporarily decorated with 1-3 and the other with 4-2.

Here are two more equivalent depictions of the same rule, coming from different choices of  1-3, 4-2 decorations. We may see the graphic beta move as

\vspace{.5cm}   \centerline{\includegraphics[width=100mm]{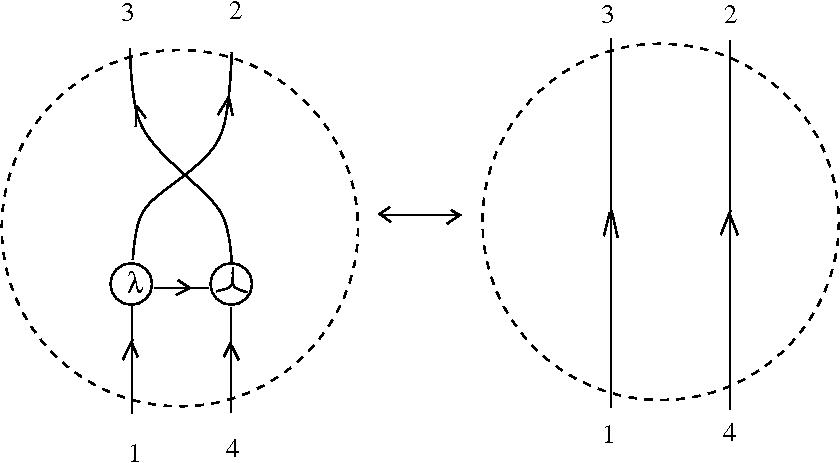}}  \vspace{.5cm}  

%    \href{http://chorasimilarity.files.wordpress.com/2012/09/betar_dif11.jpg   <img class="aligncenter  wp-image-967" title="betar_dif1" alt="" src="http://chorasimilarity.files.wordpress.com/2012/09/betar_dif11.jpg?w=300" height="278" width="506" />   

but also as

\vspace{.5cm}   \centerline{\includegraphics[width=100mm]{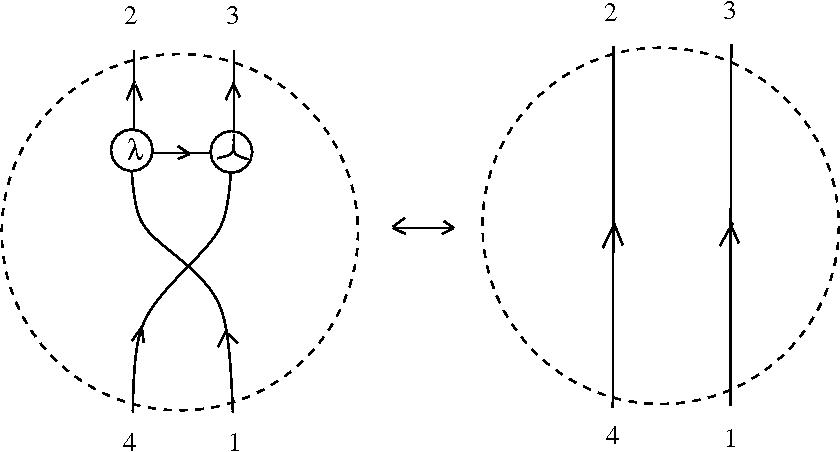}}  \vspace{.5cm}  

%    \href{http://chorasimilarity.files.wordpress.com/2012/09/betar_dif22.jpg   <img class="aligncenter  wp-image-970" title="betar_dif2" alt="" src="http://chorasimilarity.files.wordpress.com/2012/09/betar_dif22.jpg?w=300" height="277" width="517" />   

Let me then make some       notations       of the figures from the left hand sides of the previous diagrams.  Here they are: for the first figure we introduce the  (lambda calculus) "crossing macro"

\vspace{.5cm}   \centerline{\includegraphics[width=100mm]{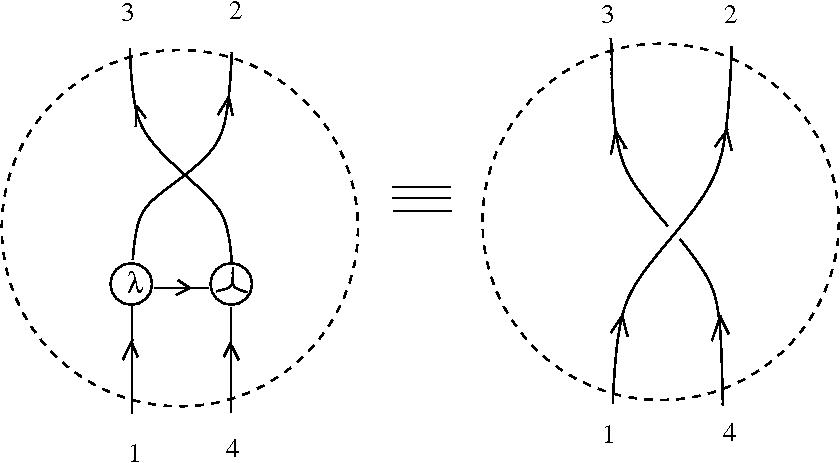}}  \vspace{.5cm}  

%    \href{http://chorasimilarity.files.wordpress.com/2012/09/betar_dif_11.jpg   <img class="aligncenter  wp-image-971" title="betar_dif_11" alt="" src="http://chorasimilarity.files.wordpress.com/2012/09/betar_dif_11.jpg?w=300" height="278" width="506" />   

and for the second figure the  other (lambda calculus) "crossing macro":

\vspace{.5cm}   \centerline{\includegraphics[width=100mm]{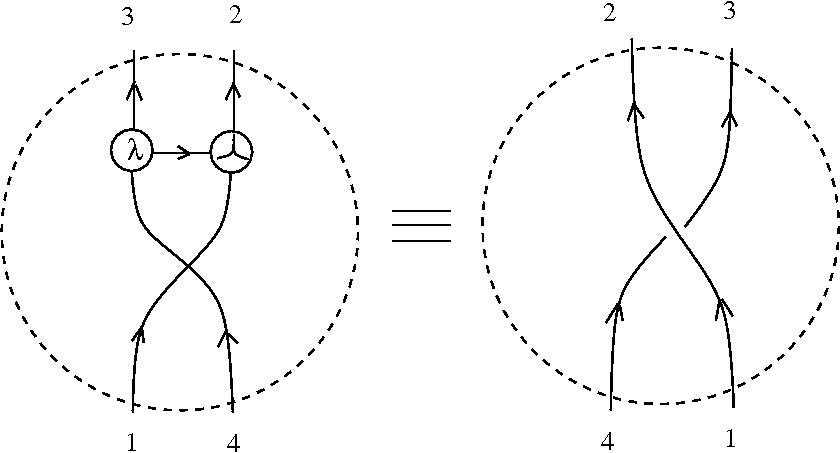}}  \vspace{.5cm}  

%    \href{http://chorasimilarity.files.wordpress.com/2012/09/betar_dif212.jpg   <img class="aligncenter  wp-image-1126" title="betar_dif21" alt="" src="http://chorasimilarity.files.wordpress.com/2012/09/betar_dif212.jpg?w=300" height="272" width="507" />   

With these notations the graphic beta move may be see as this       braiding operation

\includegraphics[width=100mm]{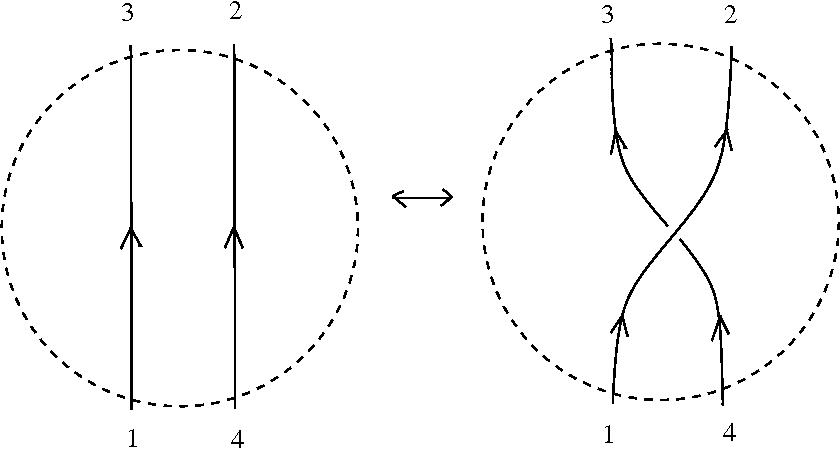}  

%    \href{http://chorasimilarity.files.wordpress.com/2012/09/betar_dif12.jpg   <img class="aligncenter  wp-image-973" title="betar_dif12" alt="" src="http://chorasimilarity.files.wordpress.com/2012/09/betar_dif12.jpg?w=300" height="259" width="486" />   

or as this other braiding operation

\vspace{.5cm}   \centerline{\includegraphics[width=100mm]{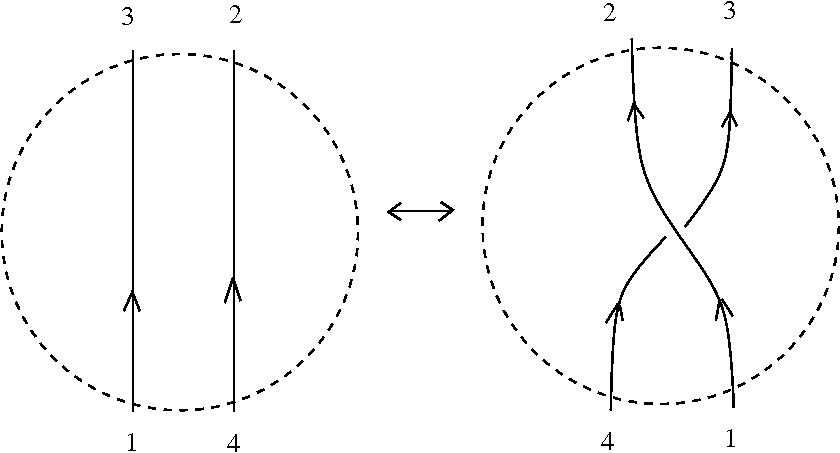}}  \vspace{.5cm}

%    \href{http://chorasimilarity.files.wordpress.com/2012/09/betar_dif222.jpg   <img class="aligncenter  wp-image-1128" title="betar_dif22" alt="" src="http://chorasimilarity.files.wordpress.com/2012/09/betar_dif222.jpg?w=300" height="277" width="521" />   

These are notations, or macros in a more technical language, which encode the performing of the graphic beta move in a particular embedding (of the graph) in the plane. As such, they look dependent of the particular embedding, but the crossings satisfy the Reidemeister moves for knot diagrams, as explained in \href{http://arxiv.org/abs/1211.1604}{arXiv:1211.1604},  \cite{graphiclambdaknots}, therefore this        notation system  is independent of the change of embedding of the graphs in the 3-dim space!

\section{Extended graphic beta move}
\label{graphicbetamoveextendedtoexplore}

       In the section \ref{emergentalgebramoves}  we defined the emergent algebra crossing macros.  It can be  noticed a very intriguing resemblance between crossings in emergent algebras (as encoded in graphic lambda) and crossings macros in graphic lambda. More specifically, here is  the encoding of a crossing in emergent algebras:

\vspace{.5cm}   \centerline{\includegraphics[width=100mm]{emerr_dif_11.jpg}}  \vspace{.5cm}

%    \href{http://chorasimilarity.wordpress.com/2012/10/31/3d-crossings-in-emergent-algebras/emerr_dif_11/" rel="attachment wp-att-1182   <img class="aligncenter size-medium wp-image-1182" alt="emerr_dif_11" src="http://chorasimilarity.files.wordpress.com/2012/10/emerr_dif_11.jpg?w=300" height="262" width="478" />   

which is very much like the crossing macro

\vspace{.5cm}   \centerline{\includegraphics[width=100mm]{betar_dif_11.jpg}}  \vspace{.5cm}  

%    \href{http://chorasimilarity.wordpress.com/2012/10/26/3d-crossings-in-graphic-lambda-calculus/betar_dif_11-2/" rel="attachment wp-att-1159   <img class="aligncenter size-medium wp-image-1159" alt="betar_dif_11" src="http://chorasimilarity.files.wordpress.com/2012/10/betar_dif_11.jpg?w=300" height="262" width="478" />

Likewise, here is the other crossing in emergent algebras:

\vspace{.5cm}   \centerline{\includegraphics[width=100mm]{emerr_dif_21.jpg}}  \vspace{.5cm} 

%    \href{http://chorasimilarity.wordpress.com/2012/10/31/3d-crossings-in-emergent-algebras/emerr_dif_21/" rel="attachment wp-att-1183   <img class="aligncenter size-medium wp-image-1183" alt="emerr_dif_21" src="http://chorasimilarity.files.wordpress.com/2012/10/emerr_dif_21.jpg?w=300" height="267" width="499" />   

and its "twin" crossing macro

\vspace{.5cm}   \centerline{\includegraphics[width=100mm]{betar_dif21.jpg}}  \vspace{.5cm}  

%    \href{http://chorasimilarity.wordpress.com/2012/10/26/3d-crossings-in-graphic-lambda-calculus/betar_dif21-3/" rel="attachment wp-att-1160   <img class="aligncenter size-medium wp-image-1160" alt="betar_dif21" src="http://chorasimilarity.files.wordpress.com/2012/10/betar_dif21.jpg?w=300" height="247" width="461" />   

These crossings are related in graphic lambda calculus through an idea introduced in  \href{http://arxiv.org/abs/1205.0139}{arXiv:1205.0139}, \cite{lambdascale}.

Here is the explanation. In $ \lambda$-Scale there are two operations, namely the $ \lambda$ abstraction and a "$\varepsilon$" operation, one for every coefficient $\varepsilon$ in a commutative group.  The application operation from lambda calculus comes as a composite of these two fundamental operations. Moreover, the $ \varepsilon$ operation IS NOT THE SAME AS the operation $\displaystyle \bar{\varepsilon}$ from emergent algebras. The emergent algebra operation $\displaystyle \bar{\varepsilon}$ is also defined as a composite of the two fundamental operations of $ \lambda$-Scale (see the paper, definition 2.2).

Later, in graphic lambda calculus, I renounced at the $ \varepsilon$ operation of $ \lambda$-Scale as a primitive operation, replacing it with the operation $\displaystyle \bar{\varepsilon}$ from emergent algebras. To be clear, I used implicitly proposition 3.2 from the mentioned paper, which gives a description of the fundamental $ \varepsilon$ operation of $ \lambda$-Scale calculus in terms of the lambda  abstraction operation, the application operation  and the emergent algebra operation $\displaystyle \bar{\varepsilon}$.

\paragraph{The extended beta move pattern.} I shall define now a macro in graphic lambda calculus (inspired by the said proposition 3.2) which corresponds to the fundamental $ \varepsilon$ operation from $ \lambda$-Scale calculus. Here is it:

\vspace{.5cm}   \centerline{\includegraphics[width=100mm]{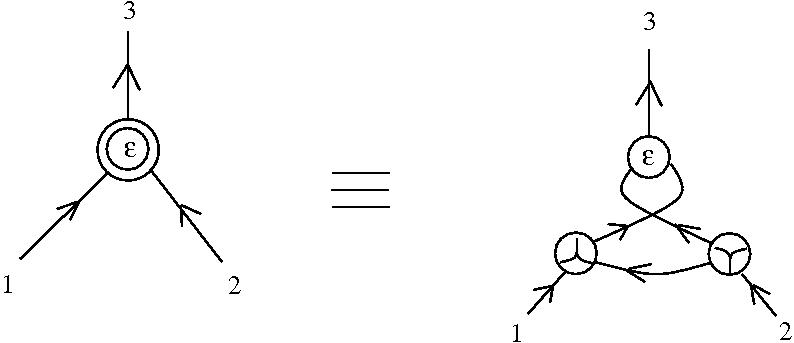}}  \vspace{.5cm}  

%    \href{http://chorasimilarity.wordpress.com/2012/12/07/lambda-scale-in-graphic-lambda-calculus-and-diagram-crossings/fundam_epsilon_1/" rel="attachment wp-att-1450   <img class="aligncenter size-medium wp-image-1450" alt="fundam_epsilon_1" src="http://chorasimilarity.files.wordpress.com/2012/12/fundam_epsilon_1.jpg?w=300" height="203" width="474" />   

We use this for defining the extended beta move pattern:

\vspace{.5cm}   \centerline{\includegraphics[width=120mm]{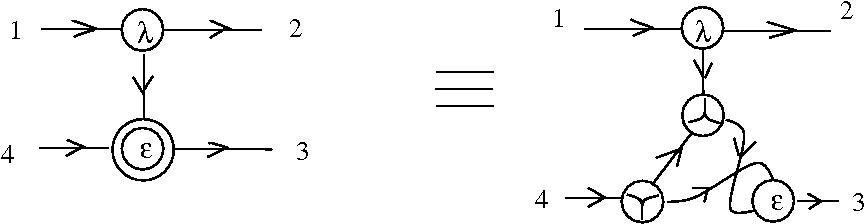}}  \vspace{.5cm}

%%%%%%%%%%%%%%%%%%%%%%%%%%%%%%%%  AICI AM AJUNS %%%%%%%%%%%%%%%%%%%%%%%%%%%%%%%%%%%%%%%%%%%%%%%%%%%%%%%%%%5

%   \href{http://chorasimilarity.wordpress.com/2012/12/07/lambda-scale-in-graphic-lambda-calculus-and-diagram-crossings/fundam_epsilon_3/" rel="attachment wp-att-1451   <img class="aligncenter size-medium wp-image-1451" alt="fundam_epsilon_3" src="http://chorasimilarity.files.wordpress.com/2012/12/fundam_epsilon_3.jpg?w=300" height="128" width="500" />   

Let's prove that the extended beta move pattern is more general than the beta move pattern and also than the pattern which appears in the emergent algebra crossing macros. 

\paragraph{The extended beta move pattern can be transformed into the beta move pattern.} Indeed, let  us consider the particular case $ \varepsilon = 1$. By using the (ext2) move and then local pruning we get:

\vspace{.5cm}   \centerline{\includegraphics[width=135mm]{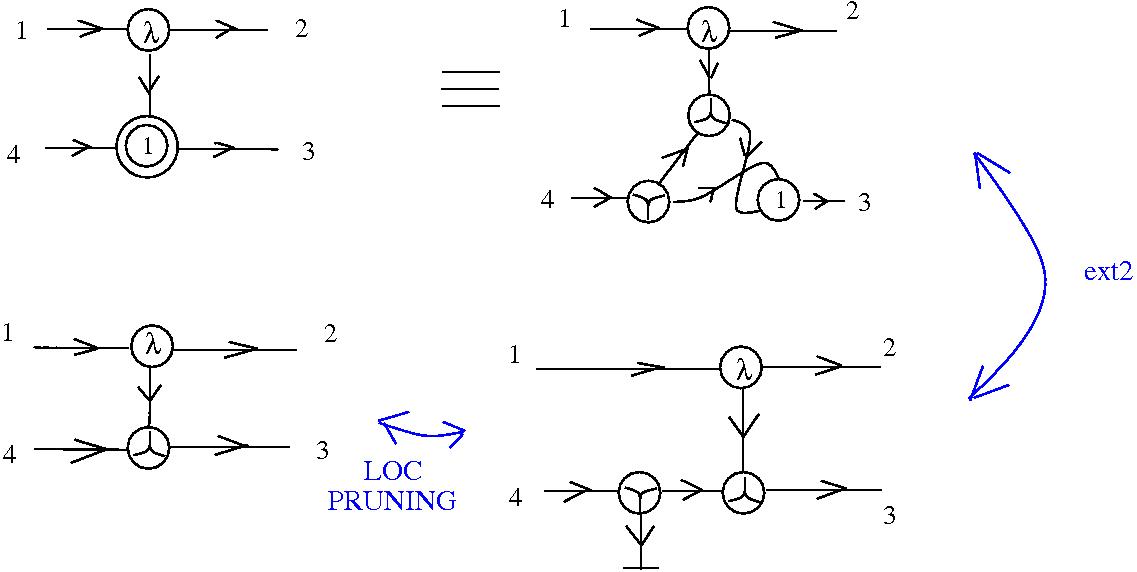}}  \vspace{.5cm} 

%    \href{http://chorasimilarity.wordpress.com/2012/12/07/lambda-scale-in-graphic-lambda-calculus-and-diagram-crossings/fundam_epsilon_4/" rel="attachment wp-att-1452   <img class="aligncenter size-medium wp-image-1452" alt="fundam_epsilon_4" src="http://chorasimilarity.files.wordpress.com/2012/12/fundam_epsilon_4.jpg?w=300" height="289" width="578" />   

\paragraph{The extended beta move pattern can be transformed into the emergent algebra crossing pattern.} I shall now apply the graphic beta move for an arbitrary $ \varepsilon$, like this:

\vspace{.5cm}   \centerline{\includegraphics[width=135mm]{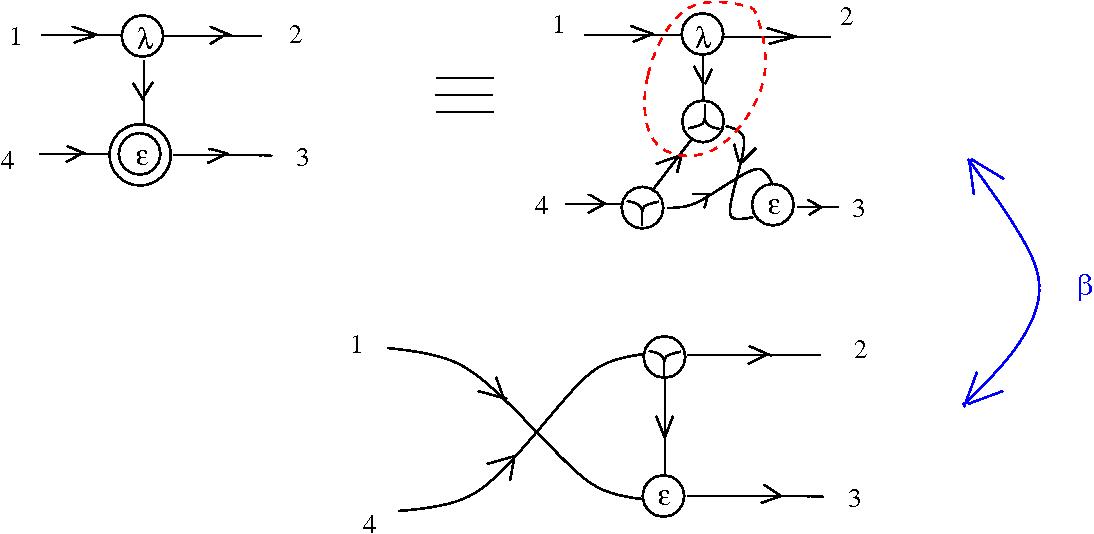}}  \vspace{.5cm}

%    \href{http://chorasimilarity.wordpress.com/2012/12/07/lambda-scale-in-graphic-lambda-calculus-and-diagram-crossings/fundam_epsilon_5/" rel="attachment wp-att-1453   <img class="aligncenter size-medium wp-image-1453" alt="fundam_epsilon_5" src="http://chorasimilarity.files.wordpress.com/2012/12/fundam_epsilon_5.jpg?w=300" height="290" width="597" />   

which is very nice, because we obtain the crossing macro of emergent algebras.

\paragraph{The extended graphic beta move.} Inspired by the introduction of the extended beta move pattern, let's define now a generalization of the graphic beta move. 

\begin{enumerate}
\item[-] The extended graphic beta move. For any $\varepsilon \in \Gamma$ we define the move "ext $\beta(\varepsilon)$" by the following figure: 
\end{enumerate}

\vspace{.5cm}   \centerline{\includegraphics[width=135mm]{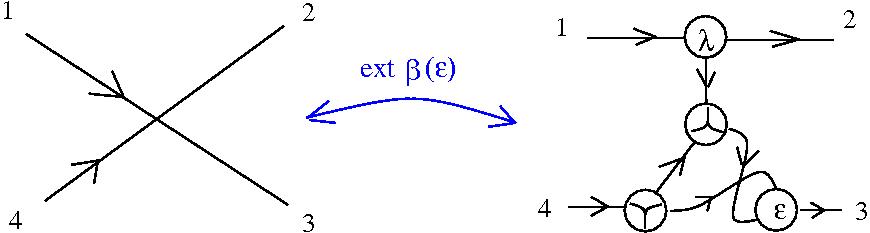}}  \vspace{.5cm}

%    \href{http://chorasimilarity.wordpress.com/2012/12/09/graphic-beta-move-extended-to-explore/fundam_epsilon_8/" rel="attachment wp-att-1462   <img class="aligncenter size-medium wp-image-1462" alt="fundam_epsilon_8" src="http://chorasimilarity.files.wordpress.com/2012/12/fundam_epsilon_8.jpg?w=300" height="163" width="616" />   

 Who could have guessed such a move, involving four gates, without thinking about diagram crossings AND about lambda calculus? Now, after expressing my enthusiasm, let's see what this move can do.

The extended graphic beta move implies the graphic beta move.  Indeed, if we take $\varepsilon =1$ then we use a previous reasoning to transform the extended beta move pattern into the graphic beta move pattern. This gives the graphic beta.

\paragraph{Dual graphic beta move.} The extended graphic beta move is equivalent with the pair of moves:

\vspace{.5cm}   \centerline{\includegraphics[width=100mm]{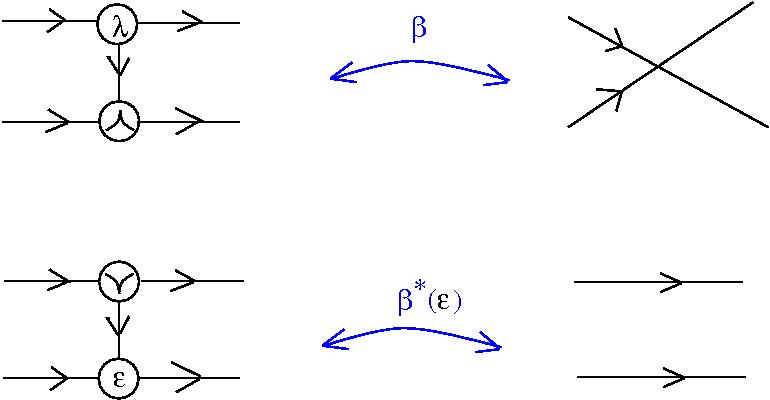}}  \vspace{.5cm}

The second move is called "the dual graphic beta move" or "$\displaystyle \beta^{*}(\varepsilon)$ move". 

Let us prove the equivalence. The extended graphic beta move implies the graphic beta move. Also, we have seen that it implies as well the $\displaystyle \beta^{*}(\varepsilon)$ move, because we can use the graphic beta move, seen as ext $\beta(1)$ move, to transform the pattern from the LHS of the ext $\beta(\varepsilon)$ move into the pattern from the LHS of the $\displaystyle \beta^{*}(\varepsilon)$ move. 

Now we have to prove that the pair (graphic beta move, $\displaystyle \beta^{*}(\varepsilon)$ move) implies the ext $\beta(\varepsilon)$. 

 Recall that in terms of the crossing macro involving only lambda calculus gates, the graphic beta move looks like this:

\vspace{.5cm}   \centerline{\includegraphics[width=100mm]{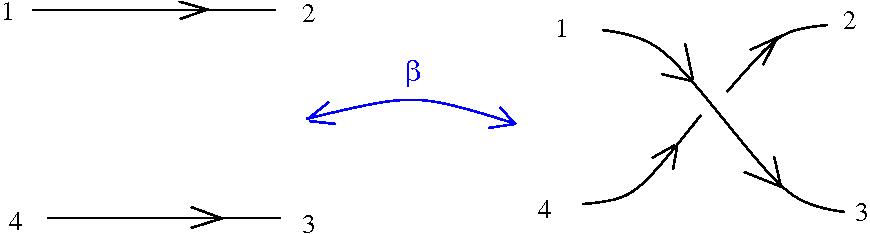}}  \vspace{.5cm} 

%    \href{http://chorasimilarity.wordpress.com/2012/12/09/graphic-beta-move-extended-to-explore/fundam_epsilon_9/" rel="attachment wp-att-1463   <img class="aligncenter size-medium wp-image-1463" alt="fundam_epsilon_9" src="http://chorasimilarity.files.wordpress.com/2012/12/fundam_epsilon_9.jpg?w=300" height="127" width="478" />   

The $\displaystyle \beta^{*}(\varepsilon)$ move, expressed via the emergent algebra crossings, has this appearance:

\vspace{.5cm}   \centerline{\includegraphics[width=100mm]{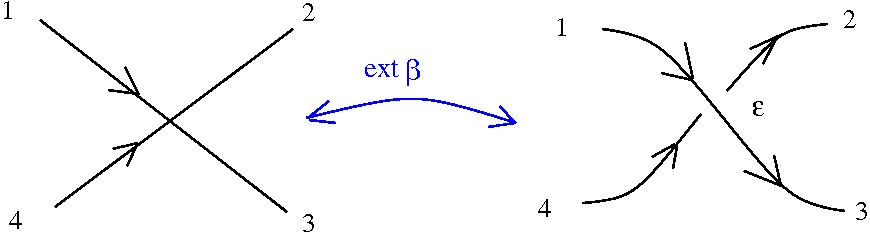}}  \vspace{.5cm}  

%    \href{http://chorasimilarity.wordpress.com/2012/12/09/graphic-beta-move-extended-to-explore/fundam_epsilon_10/" rel="attachment wp-att-1464   <img class="aligncenter size-medium wp-image-1464" alt="fundam_epsilon_10" src="http://chorasimilarity.files.wordpress.com/2012/12/fundam_epsilon_10.jpg?w=300" height="148" width="560" />   

Now, the ext $\beta(\varepsilon)$ move, expressed with crossing macros,   looks like this:

\vspace{.5cm}   \centerline{\includegraphics[width=100mm]{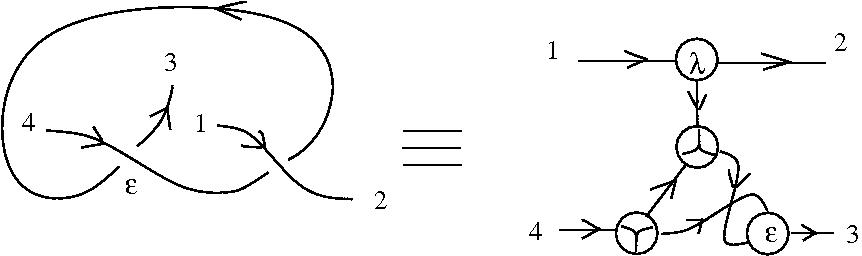}}  \vspace{.5cm} 

%    \href{http://chorasimilarity.wordpress.com/2012/12/09/graphic-beta-move-extended-to-explore/fundam_epsilon_11/" rel="attachment wp-att-1465   <img class="aligncenter size-medium wp-image-1465" alt="fundam_epsilon_11" src="http://chorasimilarity.files.wordpress.com/2012/12/fundam_epsilon_11.jpg?w=300" height="165" width="561" />   

Now we see that by applying first a $\displaystyle \beta^{*}(\varepsilon)$ move and then an usual graphic beta move, we recover the extended beta move.

\paragraph{A duality.} The equivalence between the extended graphic beta move and the pair (graphic beta move, $\displaystyle \beta^{*}(\varepsilon)$ move) suggests the following duality:

\vspace{.5cm}   \centerline{\includegraphics[width=80mm]{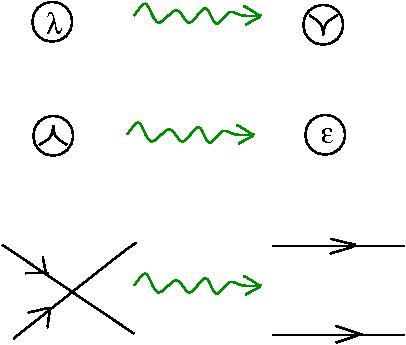}}  \vspace{.5cm}

% \href{http://chorasimilarity.wordpress.com/2012/12/17/dual-of-the-extended-beta-move-termination-as-implosion-or-blow-out/correspondence_1/" rel="attachment wp-att-1583"><img class="aligncenter size-full wp-image-1583" alt="correspondence_1" src="http://chorasimilarity.files.wordpress.com/2012/12/correspondence_1.jpg" width="406" height="345" /> 

The last correspondence should be understood like it appears in the left hand sides of the "dual" moves.

\section{Dual of the extended beta move. Termination as implosion or blow-out}
\label{dualofextendedbetamove}

If we take seriously the duality from the last section, then the dual of the extended beta move should be this:

\vspace{.5cm}   \centerline{\includegraphics[width=135mm]{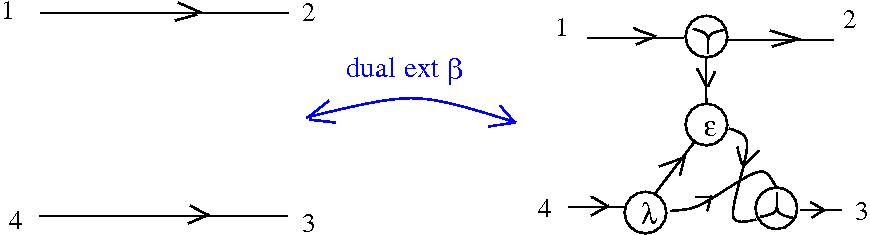}}  \vspace{.5cm}

% \href{http://chorasimilarity.wordpress.com/2012/12/17/dual-of-the-extended-beta-move-termination-as-implosion-or-blow-out/fundam_epsilon_8_dual/" rel="attachment wp-att-1584"><img class="aligncenter size-full wp-image-1584" alt="fundam_epsilon_8_dual" src="http://chorasimilarity.files.wordpress.com/2012/12/fundam_epsilon_8_dual.jpg" width="595" height="160" /> \href{

The problem is to find out if the dual of the extended beta move can be deduced from the other moves of graphic lambda calculus. If we look at the graph from the right hand side, then we see we can apply a beta move, like this:

\vspace{.5cm}   \centerline{\includegraphics[width=135mm]{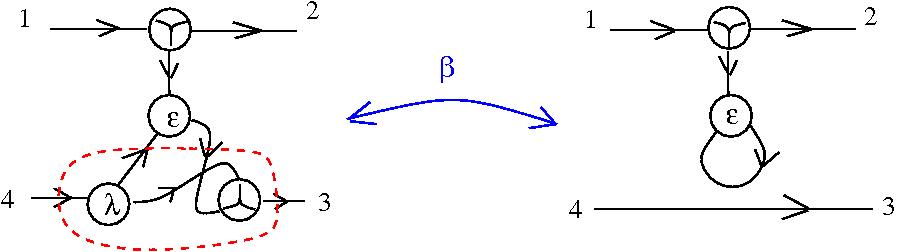}}  \vspace{.5cm}

% \href{http://chorasimilarity.wordpress.com/2012/12/17/dual-of-the-extended-beta-move-termination-as-implosion-or-blow-out/fundam_epsilon_8_dual1/" rel="attachment wp-att-1585"><img class="aligncenter size-full wp-image-1585" alt="fundam_epsilon_8_dual1" src="http://chorasimilarity.files.wordpress.com/2012/12/fundam_epsilon_8_dual1.jpg" width="595" height="167" /> \href{

Therefore the dual of the extended  extended beta move is equivalent with this curious looking move:

\vspace{.5cm}   \centerline{\includegraphics[width=135mm]{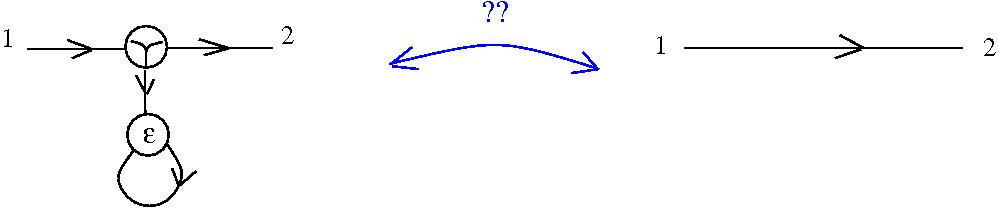}}  \vspace{.5cm}

% \href{http://chorasimilarity.wordpress.com/2012/12/17/dual-of-the-extended-beta-move-termination-as-implosion-or-blow-out/fundam_epsilon_dual2/" rel="attachment wp-att-1586"><img class="aligncenter size-full wp-image-1586" alt="fundam_epsilon_dual2" src="http://chorasimilarity.files.wordpress.com/2012/12/fundam_epsilon_dual2.jpg" width="595" height="124" /> \href{

This move looks like a  pruning move  but it is not on the list. We reason like this: by using the extended beta move and the Reidemeister 2 move (for emergent algebras) we arrive at

\vspace{.5cm}   \centerline{\includegraphics[width=110mm]{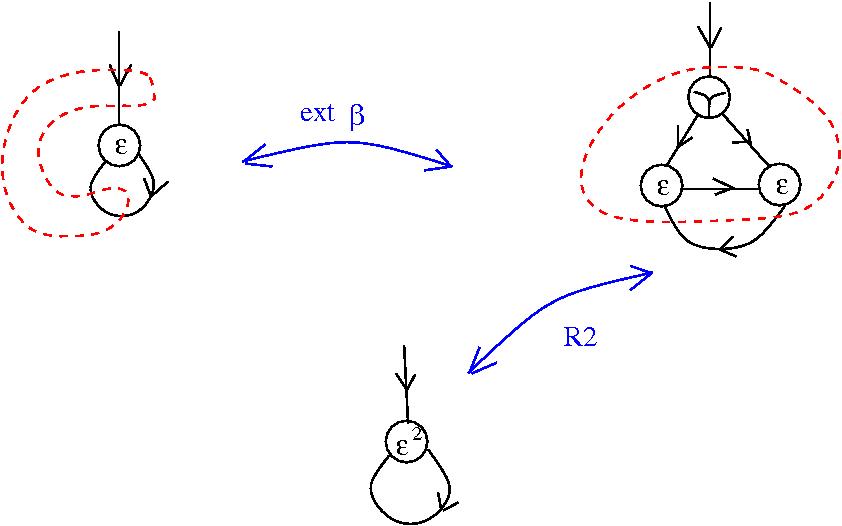}}  \vspace{.5cm}

% \href{http://chorasimilarity.wordpress.com/2012/12/17/dual-of-the-extended-beta-move-termination-as-implosion-or-blow-out/fundam_epsilon_dual3/" rel="attachment wp-att-1588"><img class="aligncenter size-full wp-image-1588" alt="fundam_epsilon_dual3" src="http://chorasimilarity.files.wordpress.com/2012/12/fundam_epsilon_dual3.jpg" width="595" height="371" /> \href{

We may repeat indefinitely this succesion of moves, or, alternatively, we may use the extended beta move with an arbitrary $ \mu$  from the group $ \Gamma$. The intuition behind this is that the gate $ \bar{\varepsilon}$ is a dilation gate, which has an input coming from the fan-out gate and the other connected to the output of the gate, therefore, by circulating along this loop, we apply an endless number of dilations of coefficient $ \varepsilon$. At the limit (if such a concept makes sense at this level of generality), either the things blow to infinity (if $ \varepsilon \geq 1$) or they implode to the value given by the fan-out gate(if $ \varepsilon \leq 1$), or they circulate forever in a loop (if $ \varepsilon = 1$) . In all cases the graph with only one input and no outputs, which we see in the previous figure, behaves exactly like the termination gate!  

Therefore, we may as well eliminate the termination gate and replace it in all the formalism by the gate from the previous figure.

\vspace{.5cm}   \centerline{\includegraphics[width=80mm]{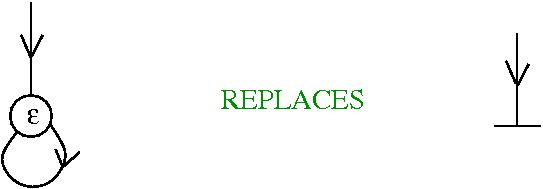}}  \vspace{.5cm}

% \href{http://chorasimilarity.wordpress.com/2012/12/17/dual-of-the-extended-beta-move-termination-as-implosion-or-blow-out/fundam_epsilon_dual4/" rel="attachment wp-att-1589"><img class="aligncenter size-full wp-image-1589" alt="fundam_epsilon_dual4" src="http://chorasimilarity.files.wordpress.com/2012/12/fundam_epsilon_dual4.jpg" width="543" height="189" /> \href{

By this replacement the "dual extended beta move" is true, as a consequence of the usual graphic beta move. Moreover, if we do this replacement, then we shall have pure trivalent graphs in the formalism because the ugly univalent termination gate is replaced by a trivalent graph!

On the other side, let us look again at the mystery move. 
The graph from the left hand side can be expressed by using the (emergent algebra) crossing macros:

\vspace{.5cm}   \centerline{\includegraphics[width=110mm]{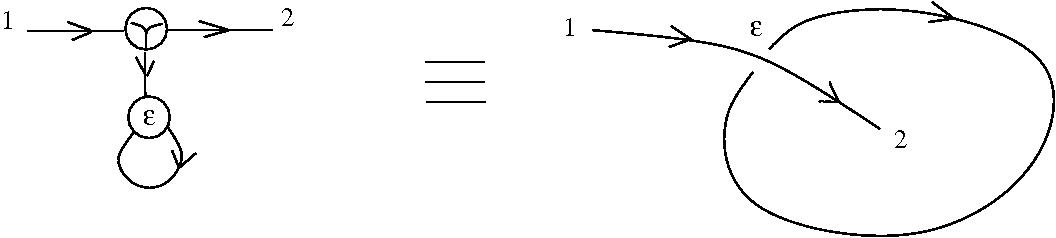}}  \vspace{.5cm} 

% \href{http://chorasimilarity.wordpress.com/2012/12/18/second-thoughts-on-the-dual-of-the-extended-beta-move/fundam_epsilon_dual5/" rel="attachment wp-att-1598"><img class="aligncenter size-full wp-image-1598" alt="fundam_epsilon_dual5" src="http://chorasimilarity.files.wordpress.com/2012/12/fundam_epsilon_dual5.jpg" width="595" height="134" /> \href{

This form of the graph makes obvious where to apply the extended beta move. Let's see what we get.

\vspace{.5cm}   \centerline{\includegraphics[width=135mm]{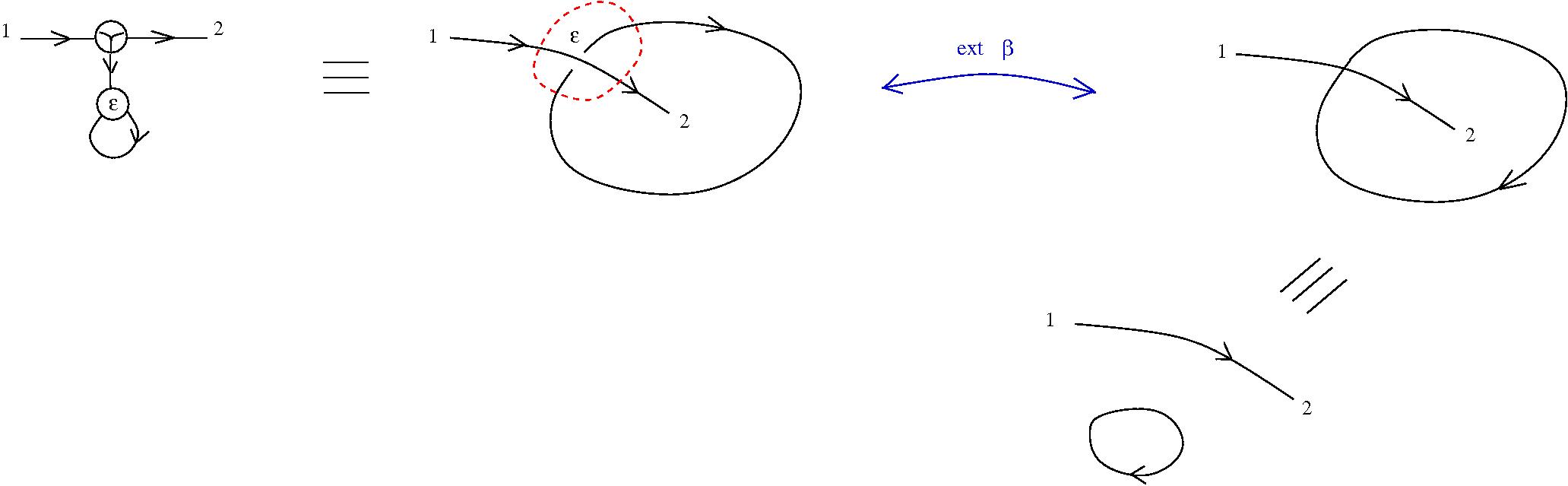}}  \vspace{.5cm} 

% \href{http://chorasimilarity.wordpress.com/2012/12/18/second-thoughts-on-the-dual-of-the-extended-beta-move/fundam_epsilon_dual6/" rel="attachment wp-att-1599"><img class="aligncenter size-full wp-image-1599" alt="fundam_epsilon_dual6" src="http://chorasimilarity.files.wordpress.com/2012/12/fundam_epsilon_dual6.jpg" width="719" height="222" /> \href{

Therefore it looks like the mystery move is a combination of the extended beta move and elimination of loops. There is no need to replace the termination gate by a trivalent graph, as suggested  in the previous post, although that is an idea worthy of further exploration.

\section{Dual of the graphic beta move implies some Reidemeister moves}
\label{dualofthegraphicbetaimpliesreidemeister}

Here are some particular applications of the dual of the graphic beta move. The first is this:

\vspace{.5cm}   \centerline{\includegraphics[width=100mm]{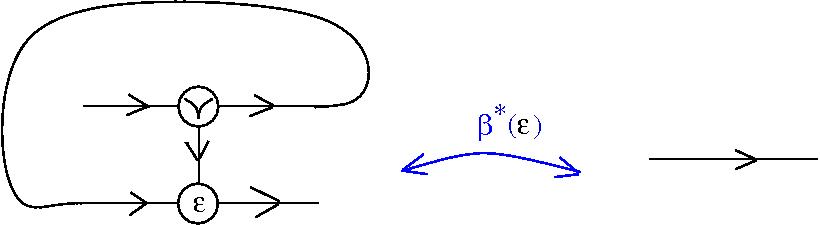}}  \vspace{.5cm} 

% \href{http://chorasimilarity.wordpress.com/2012/12/30/dual-of-the-graphic-beta-move-implies-some-reidemeister-moves/beta_beta_star_1/" rel="attachment wp-att-1688"><img class="aligncenter size-medium wp-image-1688" alt="beta_beta_star_1" src="http://chorasimilarity.files.wordpress.com/2012/12/beta_beta_star_1.jpg?w=300" width="585" height="159" /> \href{

But this is equivalent with the emergent algebra move R1a, via a  CO-COMM move.  Likewise, another application of of the dual of the graphic beta move is this:

\vspace{.5cm}   \centerline{\includegraphics[width=100mm]{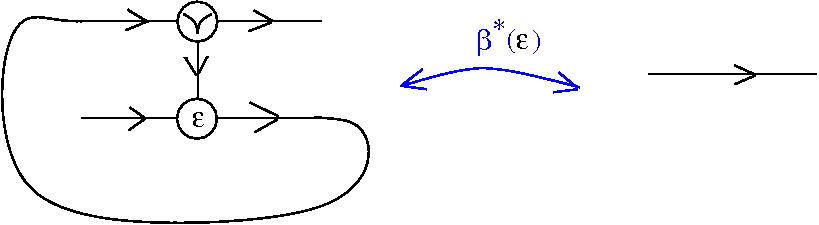}}  \vspace{.5cm}

% \href{http://chorasimilarity.wordpress.com/2012/12/30/dual-of-the-graphic-beta-move-implies-some-reidemeister-moves/beta_beta_star_2/" rel="attachment wp-att-1689"><img class="aligncenter size-medium wp-image-1689" alt="beta_beta_star_2" src="http://chorasimilarity.files.wordpress.com/2012/12/beta_beta_star_2.jpg?w=300" width="632" height="173" /> \href{

which is the same as   the  emergent algebra move R2a.

The third application is this one:

\vspace{.5cm}   \centerline{\includegraphics[width=100mm]{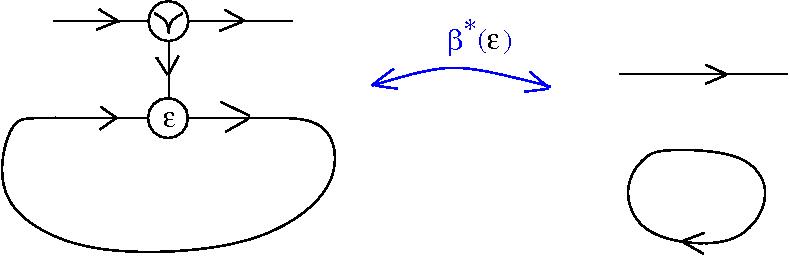}}  \vspace{.5cm}

% \href{http://chorasimilarity.wordpress.com/2012/12/30/dual-of-the-graphic-beta-move-implies-some-reidemeister-moves/beta_beta_star_3/" rel="attachment wp-att-1698"><img class="aligncenter size-full wp-image-1698" alt="beta_beta_star_3" src="http://chorasimilarity.files.wordpress.com/2012/12/beta_beta_star_31.jpg" width="595" height="192" /> \href{

which appeared  in the last part of the section \ref{dualofextendedbetamove}.

Finally, there is also this:

\vspace{.5cm}   \centerline{\includegraphics[width=100mm]{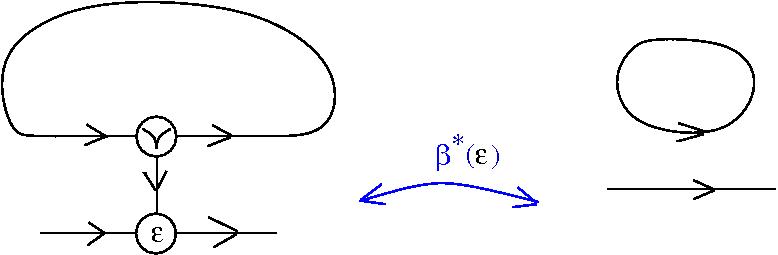}}  \vspace{.5cm}

% \href{http://chorasimilarity.wordpress.com/2012/12/30/dual-of-the-graphic-beta-move-implies-some-reidemeister-moves/beta_beta_star_4/" rel="attachment wp-att-1691"><img class="aligncenter size-medium wp-image-1691" alt="beta_beta_star_4" src="http://chorasimilarity.files.wordpress.com/2012/12/beta_beta_star_4.jpg?w=300" width="601" height="196" /> \href{

which was not mentioned before. It suggests that the following graph 

\vspace{.5cm}   \centerline{\includegraphics[width=80mm]{beta_beta_star_4.jpg}}  \vspace{.5cm} 

% \href{http://chorasimilarity.wordpress.com/2012/12/30/dual-of-the-graphic-beta-move-implies-some-reidemeister-moves/beta_beta_star_5/" rel="attachment wp-att-1692"><img class="aligncenter size-full wp-image-1692" alt="beta_beta_star_5" src="http://chorasimilarity.files.wordpress.com/2012/12/beta_beta_star_5.jpg" width="151" height="164" /> \href{

behaves like the generic point in algebraic geometry.

\end{document}